\documentclass[reqno]{amsart}

\usepackage{amsthm,amsmath,amsfonts,amssymb,epsfig,graphicx,latexsym,float,color,verbatim}
\usepackage{subcaption,tabularx,booktabs,pgfplots}
\usepackage{tikz} 

\pgfplotsset{compat=newest}
\usepgfplotslibrary{fillbetween}

\newtheorem{definition}{Definition}

\newtheorem{remark}[definition]{Remark}

\newcommand{\R}{\ensuremath{\mathbb{R}}}        

\newcommand{\G}{\ensuremath{\mathcal{G}}}
\newcommand{\Nd}{\ensuremath{\mathcal{N}}}
\newcommand{\E}{\ensuremath{\mathcal{E}}}

\newcommand{\F}{\ensuremath{\mathcal{F}}}

\begin{document}

\title[Efficient Graph-based Tensile Strength Simulations of Random Fiber Structures]{Efficient Graph-based Tensile Strength Simulations of \\ Random Fiber Structures}
\author[M.\ Harmening et al.]{Marc Harmening$^{1,\star}$ \and Nicole Marheineke$^{1}$ \and Raimund Wegener$^{2}$}

\date{\today\\
$^\star$ \textit{Corresponding author}\\
$^1$ Universität Trier, Fachbereich IV -- Mathematik, D-54286 Trier, Germany; \\ email: \{harmening, marheineke\}@uni-trier.de\\
$^2$ Fraunhofer ITWM, Fraunhofer Platz 1, D-67663 Kaiserslautern, Germany; \\ email: raimund.wegener@itwm.fraunhofer.de}

\begin{abstract}
In this paper, we propose a model-simulation framework for virtual tensile strength tests of random fiber structures, as they appear in nonwoven materials. The focus is on the efficient handling with respect to the problem-inherent multi-scales and randomness. In particular, the interplay of the random microstructure and deterministic structural production-related features on the macro-scale makes classical homogenization-based approaches computationally complex and costly. In our approach we model the fiber structure to be graph-based and of truss-type, equipped with a nonlinear elastic material law. Describing the tensile strength test by a sequence of force equilibria with respect to varied boundary conditions, its embedding into a singularly perturbed dynamical system is advantageous with regard to statements about solution theory and convergence of numerical methods. A problem-tailored data reduction provides additional speed-up, Monte-Carlo simulations account for the randomness. This work serves as a proof of concept and opens the field to optimization.
\end{abstract}

\maketitle

\noindent
{\sc Keywords.} nonwoven material; virtual tensile strength test; fiber network; graph-based data reduction; singularly perturbed dynamical system \\
{\sc AMS-Classification.}  65L04; 74Hxx; 94C15

\section{Introduction}

\begin{figure}[t]
\centering
\includegraphics[height=4cm]{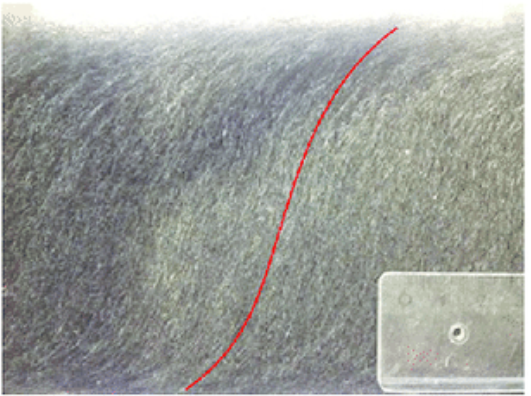}
\hspace*{2.5cm}
\includegraphics[height=4cm]{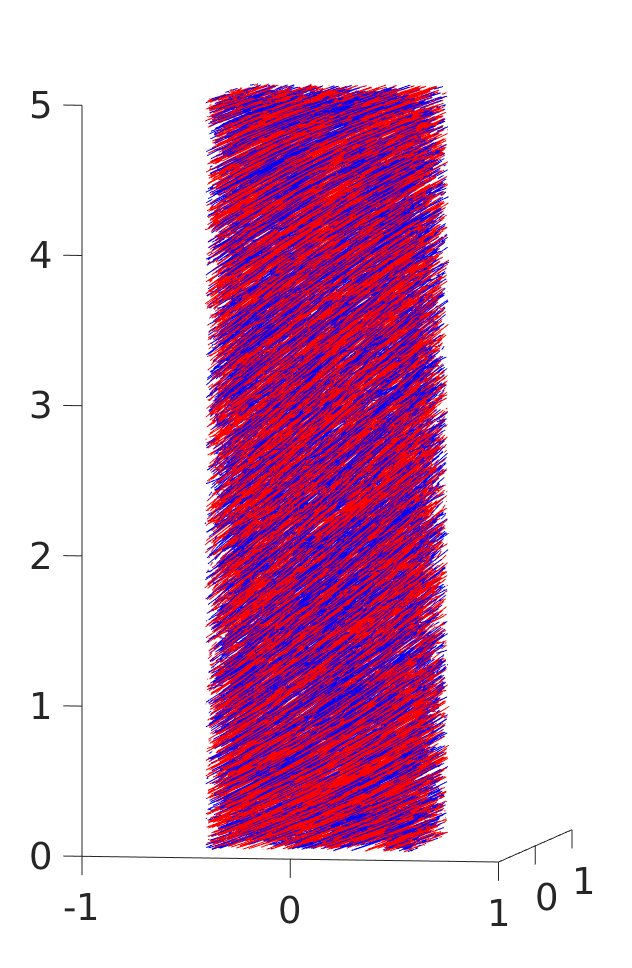}
\caption{\emph{Left:} Airlay-produced nonwoven sample that possesses a ramp-like contour plane as structural production-related feature \cite{gramsch:p:2016}. \emph{Right:} Test material volume $\mathcal{V}$ virtually generated according to the set-up specified in Sec.~\ref{sec: numerics}, Table~\ref{tab:parameter table}.}
\label{fig:tensile strength samples}
\end{figure}

Nonwovens consist of random fiber structures. Depending on the production process and the type of entering fiber material, the resulting fabrics range from hygiene and medical products over insulation and filter materials to automotive and mattress felts. The nonwovens' properties are just as diverse as their uses, they can be, e.g., homogenous, absorbent, sound-insulating, tensile. For details on nonwoven fabrics, manufacture and applications we refer to \cite{albrecht:b:2003, das:b:2014, russel:b:2006}. For the technical textile industry, the model- and simulation-based prediction and optimization of the nonwovens' properties are desirable to avoid expensive experimental test runs. Focussing here on the mechanical property of the tensile strength, the mathematical challenge lies in the appropriate handling of the existing problem-inherent multi-scales: the microstructure is determined by a random fiber network, but often superposed by structural production-related deterministic features on a macroscopic scale, see, for example, the airlay-produced nonwoven material in Fig.~\ref{fig:tensile strength samples}. We aim for a model-simulation framework that can cope with multi-scales and randomness efficiently and enable extensive tensile strength tests for design and optimization purposes.

In literature there exists a variety of ways for generating virtual material samples based on computer tomography data and/or process parameters, such as three-dimensional volume imaging \cite{faessel:p:2005, tatsuya:p:2019, ohser:b:2009}, statistical analysis and stochastic geometry \cite{ohser:b:2000, behnam:p:2006, redenbach:p:2014, schladitz:p:2006} or stochastic lay-down models \cite{goetz:p:2007,klar:p:2009,klar:p:2012}. Detailed resolution of the individual fibers, however, is cause to high complexity and makes simulations with real-size samples impossible. Instead, the use of homogenization techniques is practice where material properties are deduced from smaller and computationally handable representative volume elements \cite{hill:p:1963}. Stochastic fiber networks and non-periodic homogenization are a recent topic of research and addressed in, e.g., \cite{briane:p:1993, lebris:p:2010, lebee:p:2013, sigmund:p:1994}. For an homogenization approach on nonwoven materials see \cite{raina:p:2014}, this article also provides a remarkable survey over nonwoven microstructure models and studies in literature. Model- and simulation-based investigations of the tensile behavior and mechanical analysis of nonwoven materials can be found in, e.g., \cite{adanur:p:1999, bais:p:1995, farukh:p:2015, raina:p:2014}. Approaches based on finite element simulations on a mesoscopic scale, such as \cite{adanur:p:1999, demirci:p:2011}, treat the nonwoven as continuum and rely on empirical insights regarding the fiber orientation to incorporate randomness. Using representative microstructure samples the behavior of the fiber structure as elastic Cosserat network is computed with rods/beams in three dimensions in \cite{kufner:p:2018, strohmeyer:phd:2018} and with trusses in two dimensions in \cite{britton:p:1983, jirsak:p:1993}. In \cite{farukh:p:2015} a combination of shells and trusses is chosen to model fibrous matrix and linked single fibers. For airlay-produced nonwovens a chain of mathematical models has been presented in \cite{gramsch:p:2016} that enables the simulation of the industrial process and the investigation of the resulting nonwoven material by virtual tensile strength tests. The models range from a highly turbulent dilute fiber suspension flow to stochastic surrogates for fiber lay-down and web formation and further to random Cosserat networks with effective material laws. The nonwoven material possesses a deterministic ramp-like contour plane that crucially affects the tensile strength. To account for this macroscopic feature, representative volume elements over the material height are chosen in which the behavior of the fibers is computed by means of a finite element element method. However, in spite of homogenization and parallelization, the virtual tensile strength simulations become very costly. They are the bottleneck that the procedure is not suitable for design and optimization.

This paper aims to overcome the computational complexity. Proceeding from the model chain of \cite{gramsch:p:2016},
we equip the random graph-based fiber structure with an elastic behavior of truss-type proposing a nonlinear material law. The virtual tensile strength test is modeled by a sequence of force equilibria with respect to varied boundary conditions, yielding a family of large-scale nonlinear systems. A friction-based regularization in form of a singular perturbation ensures existence and uniqueness of solutions and makes fast, robust computations possible. Additional speed-up is achieved by the development of a problem-tailored data reduction method which removes subgraphs not contributing to the tensile behavior. Contrary to existing works, our approach is capable of handling tensile strength simulations with large complex material samples on a macroscopic level without any further homogenization techniques. This way, structural features are directly taken into account. The randomness is treated by Monte-Carlo simulations. Our work serves as a proof of concept. With the achieved efficiency the proposed framework opens the field to optimization.

The paper is structured as follows. The virtual generation of the fiber structure topology follows the spirit of \cite{gramsch:p:2016} and is briefly described in Sec.~\ref{sec:microstructure generation}. In Sec.~\ref{sec:model} we present the model for the random elastic graph-based fiber structure of truss-type and the virtual tensile strength tests. Our problem-tailored data reduction technique is subject of Sec. \ref{sec:Data Reduction}.  Using a nonwoven material of airlay-produced type, we study the performance of the tensile strength simulations in view of accuracy and computational costs in Sec.~\ref{sec: numerics}, special focus is on the effects of regularization, data reduction and randomness.

\subsection{Virtual generation of fiber structure topology} \label{sec:microstructure generation}
The nonwoven material is the deposition image of the fibers (cf.\ Fig.~\ref{fig:tensile strength samples}). In an airlay-process with short fibers, a striking feature in the fiber structure is a ramp-like contour surface. It results from the characteristic lay-down distribution of the fibers onto a moving conveyor belt in machine direction (MD), while their lay-down in cross direction (CD) and in time can be viewed as homogeneous in the production process.  A deposited fiber of length $L$ is identified with the lay-down time $T$ and the MD-CD coordinates $(X,Y)$ of one of its end points. Consider a cubic reference material volume $\mathcal{V}_R$ over the nonwoven height $h$ with base area $w_R^2$ and let $T_R$ be the time needed to produce it, i.e.\ $\mathcal{V}_R=[-w_R/2,w_R/2]^2\times[0,h]$ with $w_R>2L$. A fiber particularly contributes to the sample $\mathcal{V}_R$, if $X-x_b(T)\in [-w_R/2,w_R/2]$  is satisfied where $x_b$ accounts for the motion of the conveyor belt. In the three-dimensional web we model the fiber as a stochastic process in terms of the curve $\boldsymbol{\eta}^{(X,Y,T)}:\mathcal{I}\to \mathcal{V}_R$, $\mathcal{I}=[0,L]$
\begin{subequations}\label{sde:Final Model}
\begin{align} \label{eq:eta}
   \mathrm{d}\boldsymbol{\eta}_s &= \mathbf{R}(\boldsymbol{\eta}_s \cdot \mathbf{e_x}+x_b(T))\cdot  \boldsymbol{\tau}_s \,\mathrm{d}s , \qquad \boldsymbol{\eta}_0=(X-x_b(T))\mathbf{e_x}+Y\mathbf{e_y}+R(X)\mathbf{e_z}\\ \nonumber
    & \quad \,\, \mathbf{R}(x) = \frac{1}{\sqrt{1 + R'(x)^2}} [\mathbf{I}+(\sqrt{1 + R'(x)^2}-1)\mathbf{e_y}\otimes \mathbf{e_y} + R'(x)(\mathbf{e_z}\otimes \mathbf{e_x}-\mathbf{e_x}\otimes\mathbf{e_z})]\\\nonumber
    & \quad \,\, R(x)=h\int_{-\infty}^x r(\bar{x})\, \mathrm{d}\bar{x}
\end{align}
with $X$ $r$-distributed, $Y\sim\mathcal{U}([-w_R/2,w_R/2])$ and $T\sim \mathcal{U}([0,T_R])$ uniformly distributed -- based on the stochastic Stratonovich differential system 
    \begin{align}
     \mathrm{d} \boldsymbol{\xi}_s &=  \boldsymbol{\tau}_s\, \mathrm{d}s \label{eq:xi}\\
     \mathrm{d}  \boldsymbol{\tau}_s &= - \frac{1}{B + 1} [\boldsymbol{\Pi}_s(B)\cdot  \nabla V (\boldsymbol{\xi}_s)\, \mathrm{d}s + A \,\boldsymbol{\Pi}_s(\sqrt{B})\circ \mathrm{d} \mathbf{w}_s]\label{eq:tau} \\ \nonumber
     &\quad \, \, \boldsymbol{\Pi}_s(x)=\mathbf{n}_{\mathbf{1},s}\otimes \mathbf{n}_{\mathbf{1},s}+x\,\mathbf{n}_{\mathbf{2},s} \otimes \mathbf{n}_{\mathbf{2},s}
\end{align}
\end{subequations}
with unit tensor $\mathbf{I}$ as well as $\boldsymbol{\zeta}_0 = \mathbf{0}$ and $\boldsymbol{\tau}_0$ uniformly distributed in the unit circle spanned by $\mathbf{e_x}$ and $\mathbf{e_y}$.
The underlying anisotropic stochastic lay-down model \eqref{eq:xi}-\eqref{eq:tau} for fiber position and orientation, $((\boldsymbol{\xi},\boldsymbol{\tau}): \mathcal{I} \to \mathbb{R}^3 \times \mathbb{S}^2)$ with unit sphere $\mathbb{S}^2 \subset \mathbb{R}^3$, originates from \cite{klar:p:2012}. It presents the path of a deposited fiber onto the $\mathbf{e_x}$-$\mathbf{e_y}$ plane (here: MD-CD plane) as image of an arc-length parameterized curve that is influenced by various parameters of the production process. Modeling the fiber tangent $\boldsymbol{\tau}$, the drift term describes the typical coiling behavior with the potential $V$, while the white noise term with the vector-valued Wiener process $(\mathbf{w}:\mathcal{I} \to \mathbb{R}^3)$ and the scalar amplitude $A$ accounts for fluctuations in the lay-down process. The parameter $B \in [0,1]$ indicates the anisotropic behavior with the local orthonormal right-handed director triad $\{\boldsymbol{\tau}, \mathbf{n_1}, \mathbf{n_2} \}$, $\mathbf{n_1} \in \mathrm{span}\{\mathbf{e_x}, \mathbf{e_y}\}$, i.e., planar lay-down for $B =0$, isotropy for $B = 1$. By introducing the fiber curve $\boldsymbol{\eta}$ in \eqref{eq:eta} we adopt the standard model to feature the ramp characteristics. The contour line $R$ of the fiber material in MD is described by means of the joint probability density function $r$ of the deposited material. A fiber end point lies on the associated contour surface and the fiber orientation is aligned to it due to the local rotation group, $\mathbf{R}(x)\in SO(3)$. This model shows very concisely the typical nestling behavior of the fibers on the contour plane in contrast to the model in \cite{gramsch:p:2016} where $\mathrm{d}\boldsymbol{\eta}_s = \mathbf{R}(X) \cdot \boldsymbol{\tau}_s \, \mathrm{d}s$ with fixed rotation is used. 
In the following we restrict our considerations to the embedded test material volume $\mathcal{V}\subset \mathcal{V}_R$ with smaller base $w^2$, $w=w_R-2L$, to exclude lateral boundary effects due to partially contained fibers. The random fiber web is consolidated by adhesive joints as a result of thermobonding. Let $\boldsymbol{\eta}_\mathrm{h}$ denote the discretized fiber, i.e., set of discrete fiber points. We model an adhesive joint $\mathbf{a}$ to be formed by two fibers $\boldsymbol{\eta}_\mathrm{h}$ and $\boldsymbol{\tilde{\eta}}_\mathrm{h}$ 
\begin{align} \label{eq:contact point minimum}
&\mathbf{a}=\frac{1}{2}(\mathbf{p}^\star +\mathbf{\tilde p}^\star)\\
&\text{if} \quad \|\mathbf{p}^\star -\mathbf{\tilde p}^\star\|_2 < \kappa, \quad  (\mathbf{p}^*,\mathbf{\tilde p}^*) = \underset{(\mathbf{p},\mathbf{\tilde p}) \in \boldsymbol{\eta}_\mathrm{h} \times \boldsymbol{\tilde \eta}_\mathrm{h}}{\operatorname{argmin}}  \| \mathbf{p}-\mathbf{\tilde p} \|_2 \nonumber
\end{align}
with contact threshold $\kappa >0$. The adhesive joint takes the place of the fiber points in contact, $\mathbf{p}^\star$ and $\mathbf{\tilde p}^\star$, in the respective fibers. Note that the minimizer might be not unique such that we use the first minimizer found for practical reasons. As the fibers lie rather straight, we assume at most one contact between each fiber pair. If more fibers are involved in a contact, the resulting adhesive joint is centered between the respective fiber points in contact.  For the identification of the adhesive joints we use a bounding box method which has the complexity $\mathcal{O}((m_\eta n_\eta)^2)$, $m_\eta$ discretization points per fiber and $n_\eta$ fibers in total.

Of course, the virtual generation of microstructures that are composed of different fiber types is straightforward possible. Then, the initial MD coordinate $X$ of $\boldsymbol{\eta}$ is distributed wrt.\ the MD lay-down density function $g_j$ of the respective fiber type and the joint probability density function becomes $r=\sum_{j} \beta_j g_j$, $\sum_j \beta_j=1$ where the weights $\beta_j\geq 0$ are determined wrt.\ the mass ratio or the number ratio of the fiber types. Different adhesive properties of the fiber types can be taken into account in the identification of adhesive joints, see \cite{gramsch:p:2016} for details.

\section{Modeling Framework}\label{sec:model}
The virtual fiber structure can be represented by a random variable $\mathcal{M} = \{ \mathcal{G}, \F, \boldsymbol{\ell}, \mathbf{x} \}$. The topology of the structure is described by a graph $\mathcal{G}(\mathcal{N},\mathcal{E})$ where the nodes $\Nd$ represent adhesive joints and fiber ends and the edges $\E$ indicate the existence of fiber connections between them. In particular, multiple fiber connections between two nodes are represented by a single edge.  The set of all fiber connections is denoted by $\F$, $| \F | \geq |\E |$. Each fiber connection is equipped with a fixed (positive) length, yielding the gobal length vector $\boldsymbol{\ell}\in \mathbb{R}_+^{|\mathcal{F}|}$.  The spatial positions of the adhesive joints and fiber ends (nodes) are given by $ \mathbf{x}\in\R^{3| \Nd |}$. To refer to the position of an individual node $\nu \in \Nd$ we subsequently write $ \mathbf{x}_{\nu} \in \R^3$. Note that we distinguish between the nodes $\mathcal{N}_{B}$ associated to the upper and lower faces (boundary) of the fiber structure  and the remaining interior nodes $\mathcal{N}_I$, i.e., $\mathcal{N}=\mathcal{N}_I \dot \cup \,\mathcal{N}_{B}$. In the following we focus on a single sample $M=\mathcal{M}(\omega)$. The randomness is accounted for by Monte-Carlo simulations.

In the experiment the tensile strength of a nonwoven structure is determined by the stress-strain relation when pulling the upper and lower faces apart. Aiming for efficient tensile strength simulations we equip the virtual graph-based structure with an elastic behavior of truss-type, neglecting angular momentum effects. We mimic the experiment by considering a sequence of force equilibria with respect to varied boundary conditions, yielding a family of large-scale nonlinear equation systems. A regularization in form of a dynamic embedding ensures existence and uniqueness of solutions and makes fast, robust computations possible. 

\subsection{Elastic graph-based fiber structure of truss-type}  \label{sec:modelling}
The elastic material behavior can be modeled by treating the fiber connections as Cosserat rods or Timoshenko beams as, e.g., done in \cite{strohmeyer:phd:2018}. However, the quality of such a model is in no way related to the performance. The angular effects are very small, but the computational effort is enormous so that additional homogenization techniques have to be used for realizing tensile strength simulations. Aiming for computations on the macroscopic scale that are fast and robust in view of Monte-Carlo simulations we pursue a different strategy. We equip the fiber structure with an elastic behavior of nonlinear truss-type. The present constellation of the nodes and, consequently, of the whole structure is then determined by the acting (traction) forces on the fiber connections.

\subsubsection{Equilibrium model} 
\begin{figure}
\centering
    \includegraphics[width=0.5\linewidth]{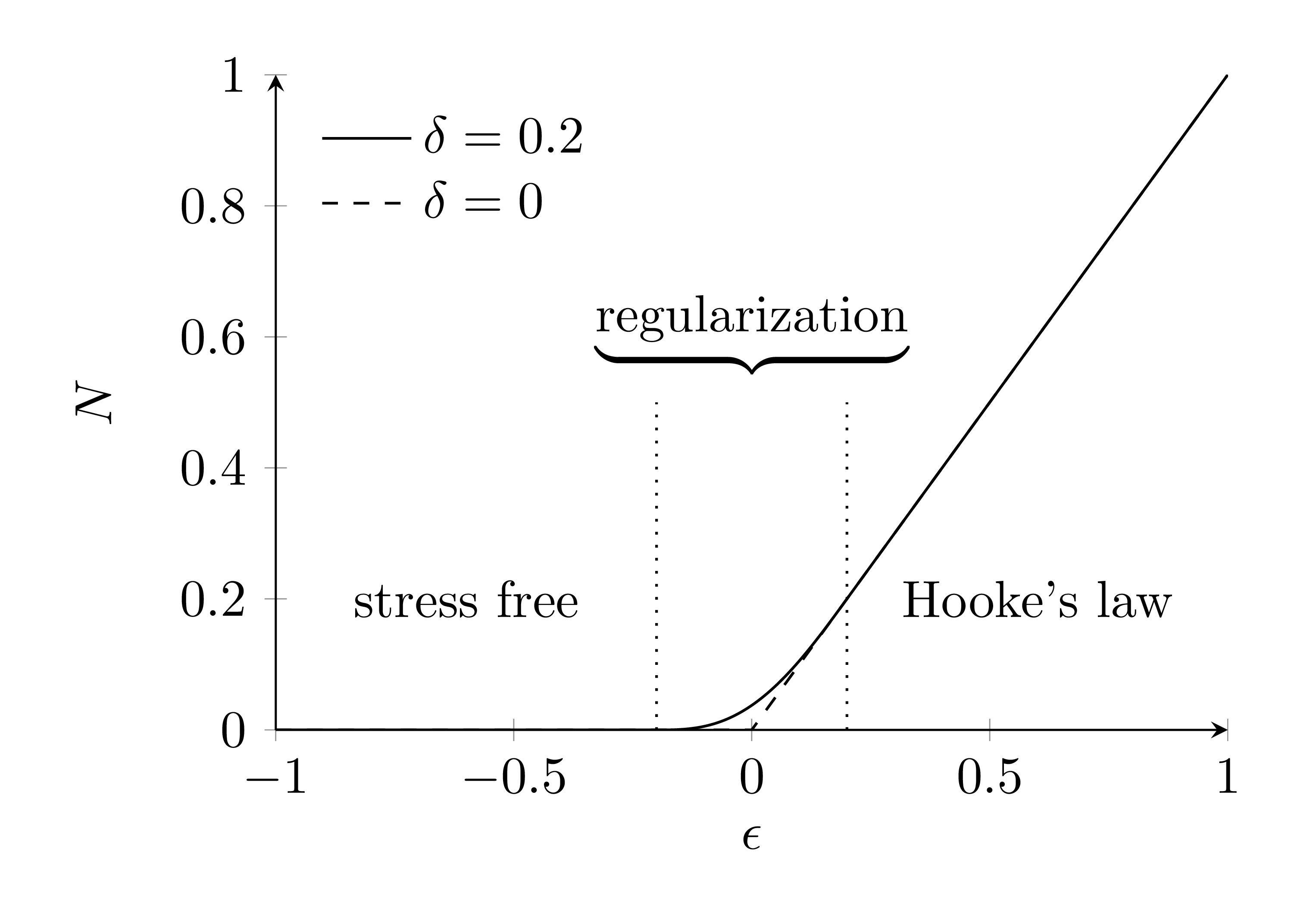}
    \caption{Material law \eqref{eq:stress-strain relation} for $EA = 1$~[N] and $\delta = 0.2$.}
    \label{fig:Fiber Law}
\end{figure}

Assume a static set-up. Given the positions of the boundary nodes $\Nd_{B}$, a force equilibrium with respect to the acting (traction) forces on the incident fiber connections is satisfied in all interior nodes $\Nd_I$, i.e.,
\begin{subequations}\label{eq:equilibrium}
\begin{align} \label{eq:fixed node position}
 \mathbf{x}_{\nu} &= \mathbf{g}_{\nu}, && \forall \nu \in \Nd_B \\ \label{eq:force equilibrium}
\sum\limits_{\mu \in \E(\nu)} \sum\limits_{r \in \F(\mu)} \mathbf{f}_{\mu, r}^{\nu}(\mathbf{x}) &= \mathbf{0}, \qquad \mathbf{f}_{\mu, r}^{\nu}(\mathbf{x}) =    
    \frac{\mathbf{t}_{\mu}^{\nu}(\mathbf{x})}{\| \mathbf{t}_{\mu}^{\nu}(\mathbf{x}) \|_2} N(\,\epsilon(\,\| \mathbf{t}_{\mu}^{\nu}(\mathbf{x}) \|_2 ,\, \ell_{r}\,)\,),  && \forall \nu \in \Nd_I   
\end{align}
where $\E(\nu)\subset E$ denotes the edges being incident to node $\nu$ and $\F(\mu)\subset \F $ the fiber connections being represented by edge $\mu$. The force $\mathbf{f}_{\mu, r}^{\nu}: \R^{3|\mathcal{N}|} \to \R^3$ in node $\nu$ that is caused by fiber connection $r$ belonging to edge $\mu$ acts in normalized edge direction with  $\mathbf{t}_{\mu}^{\nu}(\mathbf{x})= \mathbf{x}_{\Tilde \nu} - \mathbf{x}_{\nu}$ for $\mu = (\nu,\Tilde \nu)$. Its amplitude $N$ depends on the relative strain of the fiber connection with respect to its length $\ell_r$, i.e., $\epsilon: \mathbb{R}^+ \times \mathbb{R}^+ \rightarrow (-1,\infty)$, $(l,\ell) \mapsto (l-\ell)/\ell$.  We model it as 
\begin{align} \label{eq:stress-strain relation}
N(\epsilon) = \left\{
\begin{array}{cl}
0,   &-1 < \epsilon \leq -\delta\\
EA \,\epsilon \, \bigl( -\frac{1}{16} (\frac{\epsilon}{\delta})^3 + \frac{3}{8}(\frac{\epsilon}{\delta}) + \frac{1}{2} + \frac{3}{16} (\frac{\epsilon}{\delta})^{-1} \bigr),   &-\delta < \epsilon \leq \delta \\
EA \, \epsilon,  &\phantom{-}\delta < \epsilon   \\
\end{array}
\right. 
\end{align}
\end{subequations}
which reflects Hooke's law, the common linear stress-strain relation with respect to the cross-section associated Young's modulus $EA$, in the stretched state ($\epsilon >0$). In the unstretched state ($\epsilon<0$) we assume that the force required to straighten the fiber is negligibly small and hence set it to be zero. The interpolating polynomial on $[-\delta,\delta]$ with $0<\delta \ll 1$ is meant to rise the regularity such that $N \in  \mathcal{C}^2((-1,\infty), \R^+_0)$ (see Fig.~\ref{fig:Fiber Law}). 

The variables of the equilibrium model \eqref{eq:equilibrium} are the positions of the interior nodes $\mathbf{z}\in \R^{3|\mathcal{N}_I|}$. Hence we introduce $\mathbf{x}(\mathbf{z})=\mathbf{ P_I}^T \cdot \mathbf{z} + \mathbf{P_B}^T \cdot \mathbf{g}$ with the projections $\mathbf{ P_I}\in \R^{3|\mathcal{N}_I| \times 3|\mathcal{N}|}$ and $\mathbf{ P_B} \in \R^{3|\mathcal{N}_B| \times 3|\mathcal{N}|}$ onto the interior and boundary nodes, respectively. Rewriting and scaling the model, it is convenient to consider it as a dimensionless nonlinear system of the form 
\begin{align}\label{eq:NES}
\mathbf{F}(\mathbf{z})=\mathbf{0}, \qquad \mathbf{F}:\R^{3|\mathcal{N}_I|} \to \R^{3|\mathcal{N}_I|}
\end{align}  
which we refer to as
Network Equation System (NES). The function $\mathbf{F}$ inherits the regularity of $N$, in particular it and its Jacobian are Lipschitz continuous. Also related to the specific choice of the material law is that a solution of NES is in general not unique. To illustrate this we consider the following example: Let $\mathbf{z}$ be a solution of NES. We assume that there exists a node $\nu \in \Nd_I$ that is only incident to fiber connections in an unstretched state, i.e., $\epsilon(\|\mathbf{t}^\nu_{\mu}(\mathbf{x}(\mathbf{z}))\|_2,\ell_r) \in (-1,-\delta]$ for all $\mu \in \E(\nu)$ and $r \in \F(\mu)$. Since $N(\epsilon)=0$ holds for  $\epsilon \in (-1,-\delta]$, sufficient small variations in the position of $\nu$ do not exert stress on any of the incident fiber connections. Hence all these variations also solve NES. The lack of uniqueness comes from the fact that $N$ is not injective. Injectivity could be easily achieved, for example, by introducing crimp as done in \cite{gramsch:p:2016}, but such modeling is rather heuristic. We prefer the zero-phase in $N$, as it yields computational advantages concerning data reduction and sparsity of the Jacobian of $\mathbf{F}$ when dealing with the large-scale system after appropriate regularization.

\subsubsection{Quasi-static model with friction-based regularization} \label{subsec:dynamic embedding}
To model the tensile strength test we consider a family of force equilibria \eqref{eq:NES} with respect to varied boundary positions, i.e., for $t\in [0,1]$ 
\begin{align*}
    \mathbf{x}_{\nu} = \mathbf{g}_{\nu}(t),  \qquad  \mathbf{g}_{\nu}(t)&= {\mathbf{x}}^\circ_{\nu} + t \,\mathbf{d}_{\nu} && \forall \nu \in \Nd_B\\
     \mathbf{d}_\nu&=\left\{\begin{array}{l l}
     d \, \mathbf{e}_3, \,\,d>0, & \nu \in \mathcal{N}_{B_u}\\
     \mathbf{0}, & \nu \in \mathcal{N}_{B_l}.
     \end{array}\right.
\end{align*}
Proceeding from an initial constellation ${\mathbf{x}}^\circ_{\nu}$, the nodes associated to the upper face ($\nu \in \mathcal{N}_{B_u}$) are linearly shifted away with (maximal) displacement length $d>0$, while the lower face is fixed ($\nu \in \mathcal{N}_{B_l}$). This yields the family of large-scale nonlinear systems $\mathbf{F}(\mathbf{z},t)=\mathbf{0}$, $t \in [0,1]$ which we refer to as quasi-static model.  The regularity of $\mathbf{F}$ might provide quadratic convergence of Newton's method when numerically solving the quasi-static model as  a sequence of equilibrium problems for $0=t_0<t_1<...<t_M=1$ in a continuation framework and using the solution $ \mathbf{z}(t_i)$ as initial guess for $\mathbf{F}(\mathbf{z}, t_{i+1})=\mathbf{0}$. However, the lack of uniqueness might cause singular Jacobians which prevents a straightforward application of the method. 

To overcome the flaw while preserving the convergence properties we introduce a regularization. We propose a regularized (dimensionless) model variant in form of a singularly perturbed (stiff) system of ordinary differential equations with $\varepsilon\ll 1$, for $t \in (0,1]$
\begin{equation}\label{eq:dimensionless dynamic formulation}
\begin{aligned} 
\varepsilon  \, \dot{\mathbf{z}}(t) & = \mathbf{F} (\mathbf{z}(t),t), \qquad &&
\mathbf{z}(0) = \mathbf{P_I} \cdot \mathbf{x}^\circ  .
\end{aligned}
\end{equation}
The regularization term can be interpreted as friction that is experienced by the interior nodes in the fiber structure during the tensile strength test. In this context the regularization parameter stands for $ \varepsilon = \gamma \, w /((EA) \,T)$ with fiber friction coefficient $\gamma$, elasticity modulus $EA$, test volume width $w$ and final test time $T$. As the tensile strength test takes long time (the fiber structure is pulled away very slowly), $\varepsilon \ll 1$ holds true and the transition to the quasi-static framework is justified, $\varepsilon \rightarrow 0$. Due to the regularity of $\mathbf{F}$ and in particular its Lipschitz continuity in $\mathbf{z}$, the regularized model \eqref{eq:dimensionless dynamic formulation} possess a unique solution $\mathbf{z}_\varepsilon$. Moreover, $\mathbf{z}_\varepsilon(t)$ approaches to a solution of $\mathbf{F}(\mathbf{z},t)=\mathbf{0}$ for all $t \in [0,1]$ as $\varepsilon \rightarrow 0$. This might be interpreted as the friction-associated solution of the quasi-static model. For the numerical handling of \eqref{eq:dimensionless dynamic formulation} we apply an implicit integration scheme to account for the model's stiffness and solve the resulting nonlinear systems with Newton's method. In each time step we use solution contributions associated to the former time levels as initial guess.

From modeling side the consideration of the singularly perturbed system \eqref{eq:dimensionless dynamic formulation} is intuitive, resulting in a good numerical performance concerning efficiency and accuracy of the tensile strength simulations. From optimization side also other regularizations of \eqref{eq:NES} might be possible, e.g., a classical Tikhonov regularization \cite{hou:p:1997,ito:b:2014}. In that case the uniqueness is obtained by adding a (small physically artificial) force in direction of some reference node constellation.  Alternatively, other numerical approaches, such as quasi Newton methods, or minimal residual techniques that are able to handle singularities in the Jacobian, could be applied directly to NES, but they possess in general lower convergence rates \cite{hanke:p:1995, jin:p:2000}.

\begin{remark} \label{remark:model extension}
In case that the fiber structure is viewed as a representative volume element in a homogenization analysis, it might make sense to keep the lateral nodes on the respective faces in the tensile strength model such that a lateral contraction of the structure is prevented during the simulation. The lateral faces are particularly induced by the base $ [-w/2,w/2]^2$ of the initially surrounding test volume $\mathcal{V}$ with outer normals $\mathbf{n}$. Let $\Nd_L $ be the set of the lateral nodes such that $\Nd=\mathcal{N}_I \dot \cup \,\mathcal{N}_{B} \dot \cup \,\mathcal{N}_{L}$ holds. Then, \eqref{eq:equilibrium} is supplemented by the following conditions
\begin{subequations}
\begin{align}
     \sum\limits_{\mu \in \E(\nu)} \sum\limits_{r \in \F(\mu)} \mathbf{f}_{\mu, r}^{\nu}(\mathbf{x}) + \lambda_{\nu} \mathbf{n}_{\nu} &= \mathbf{0}, && \forall \nu \in \Nd_L \label{eq:lagrange multiplicator}\\
     \mathbf{x}_{\nu} \cdot \mathbf{n}_{\nu} - \frac{w}{2} &= 0 , && \forall \nu \in \Nd_L \label{eq:boundary fixation}
     \end{align}
\end{subequations}
where the force equilibrium in a lateral node $\nu$ is augmented by a force acting in the direction of the respective outer normal $n_{\nu} \in \R^3$. The force amplitude $\lambda_{\nu} \in \R$ acts as Lagrange multiplier (unknown) to the location constraint on $\nu$. The modified equilibrium model with the positions of interior and lateral nodes $\mathbf{z}\in \R^{3(|\mathcal{N}_I|+|\mathcal{N}_L|)}$ as well as the Lagrange multipliers $\boldsymbol{\lambda} \in \R^{ |\Nd_L |}$ as variables can be formulated in (dimensionless) form
\begin{equation}\label{eq:NES_2}
\begin{aligned}
\Tilde{\mathbf{F}}(\mathbf{z},\boldsymbol{\lambda}) &=\mathbf{0}, \qquad && \Tilde{\mathbf{F}}:\R^{3(|\mathcal{N}_I|+|\mathcal{N}_L|)}\times \R^{|\mathcal{N}_L|} \to \R^{3(|\mathcal{N}_I|+|\mathcal{N}_L|)}\\
\mathbf{G}\cdot\mathbf{z}&=\mathbf{0}, &&  \mathbf{G}\in \R^{|\mathcal{N}_L| \times 3(|\mathcal{N}_I|+|\mathcal{N}_L|)}
\end{aligned}
\end{equation}
where the nonlinear function $\Tilde{\mathbf{F}}$ comes from \eqref{eq:equilibrium}, \eqref{eq:lagrange multiplicator} and the linear mapping from the location constraint \eqref{eq:boundary fixation}. The matrix $\mathbf{G}$ has full row rank. Analogously to \eqref{eq:NES}, the system \eqref{eq:NES_2} possesses multiple solutions.

The corresponding regularized quasi-static model variant becomes a semi-explicit system of differential algebraic equations, for $t\in (0,1]$
\begin{equation} \label{eq:semi-explicit DAE}
\begin{aligned}
\varepsilon\, \dot{\mathbf{z}}(t)&=\Tilde{\mathbf{F}}(\mathbf{z}(t),\boldsymbol{\lambda}(t),t), \qquad && \mathbf{z}(0) = \mathbf{P_{I,L}} \cdot \mathbf{x}^\circ\\
  \mathbf{0} &= \mathbf{G} \cdot \mathbf{z}(t), && \boldsymbol{\lambda}(0)=\mathbf{0}
\end{aligned}
\end{equation}
since the regularization (adding of friction) is only reasonable in the force balances whereas the location constraints are kept unchanged as algebraic equations. Here,  $\mathbf{P_{I,L}}$ denotes the projection onto the interior and lateral node positions. The system \eqref{eq:semi-explicit DAE} is of differentiation index~2.
\end{remark}

\begin{figure}[t]
\begin{subfigure}[t]{0.39\textwidth}
    \centering
    \includegraphics[width = \linewidth]{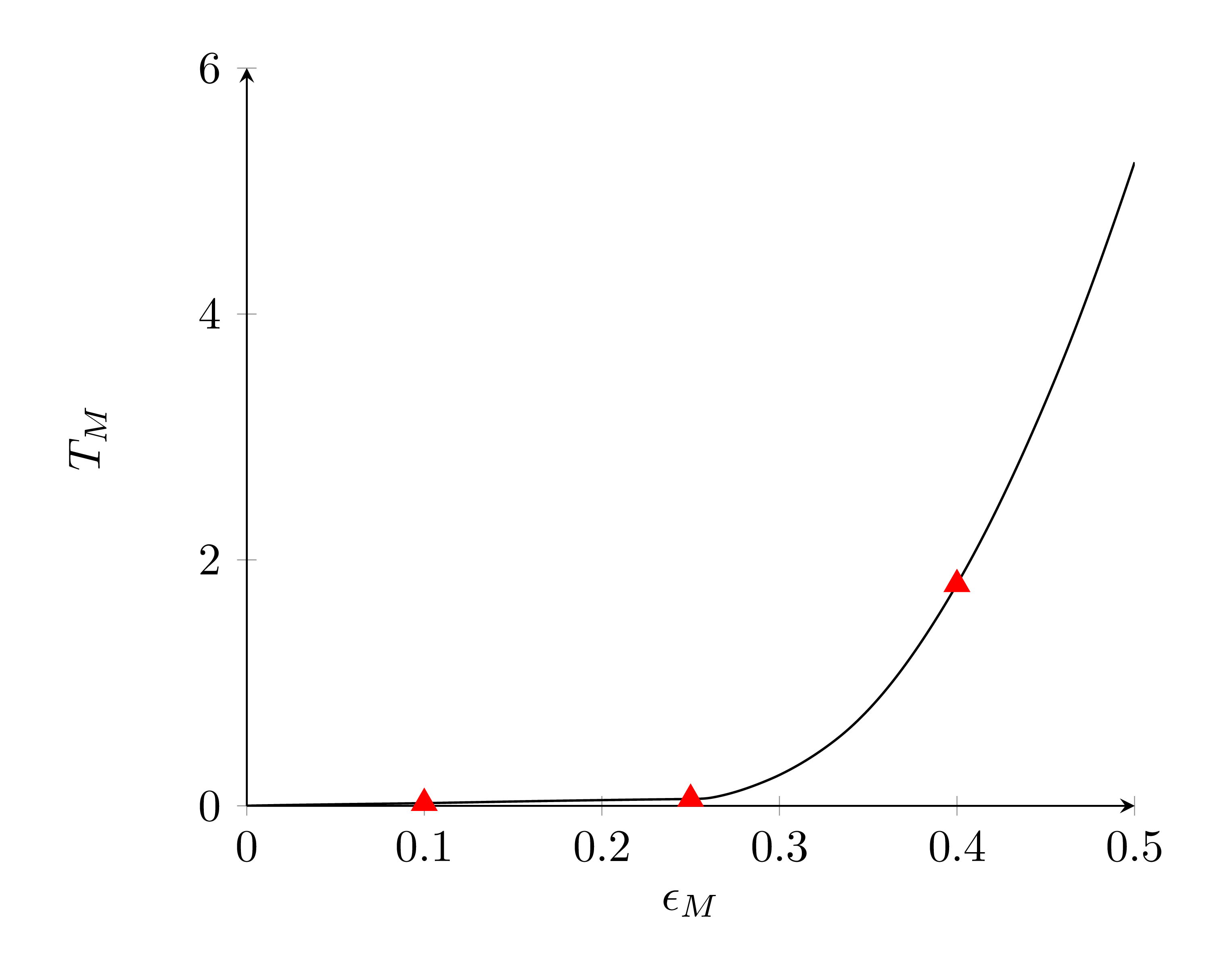}
    \caption{\tiny Stress-strain relation with tensile force $T_M$ [N].}
    \label{fig:Stress Strain Curve}
\end{subfigure}
\hfill
\begin{subfigure}[t]{0.56\textwidth}
    \centering
    \includegraphics[height=4.65cm]{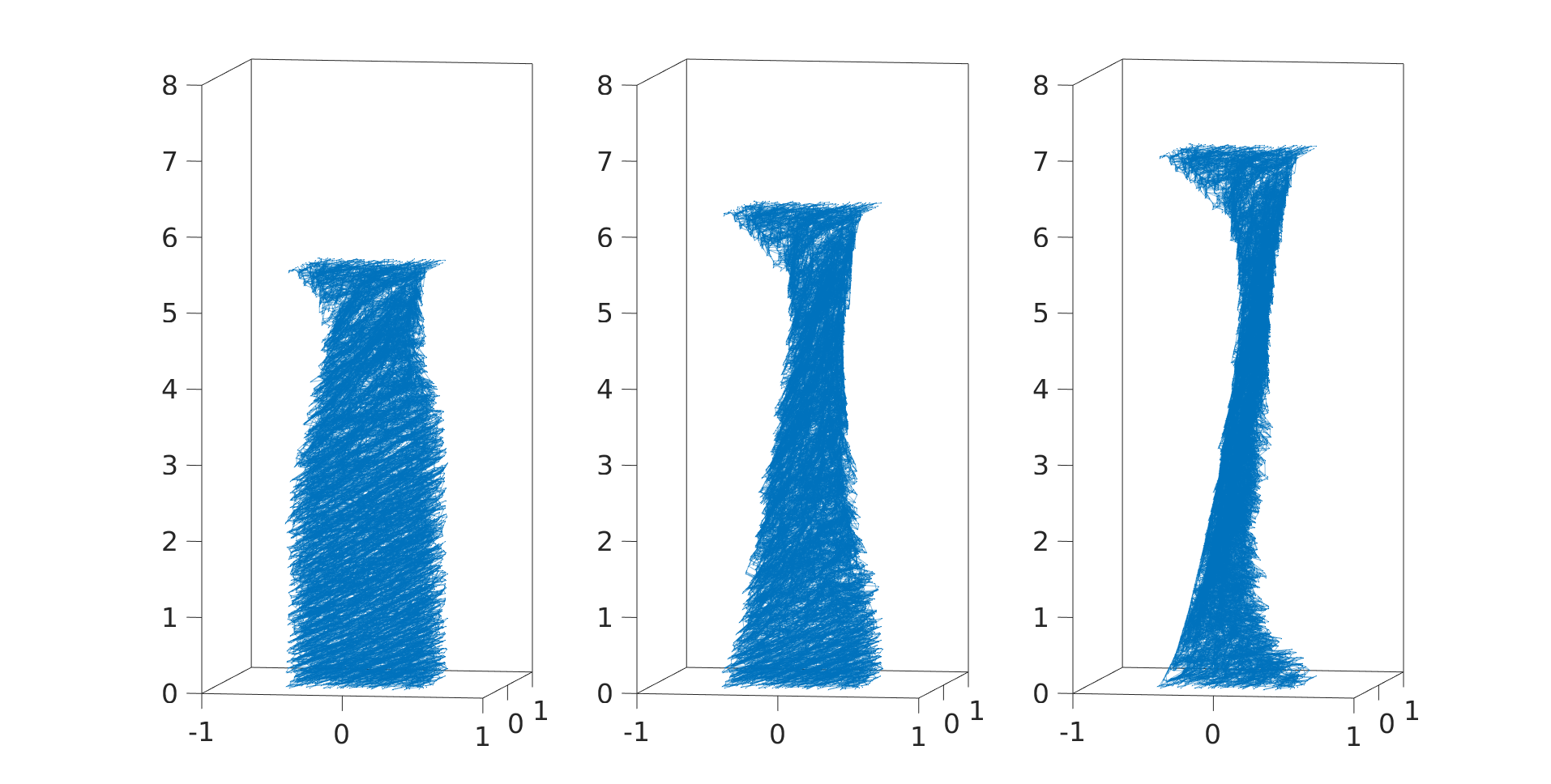}
    \caption{\tiny Behavior of virtual fiber structure during the test at $(\epsilon_M,T_M)$ \\ \phantom{.......} indicated by markers in (A). }
    \label{fig:Test volume bahvior}
\end{subfigure}
\caption{Tensile strength simulation results computed with \eqref{eq:dimensionless dynamic formulation} for fiber structure sample of Fig.~\ref{fig:tensile strength samples}.}
\label{fig:exemplary simulation}
\end{figure}

\subsection{Tensile strength behavior} \label{sec:tensile strength behavior}
Given the fiber structure sample $M=\mathcal{M}(\omega)$, its tensile strength behavior is quantified by its stress-strain relation. Let $\mathbf{e_3}$ be the shift direction in the tensile strength experiment.  By pulling away the upper face  while fixing the lower one of the microstructure, an inner tensile force is caused, which acts against the pulling. Its amplitude adds up from all force components of fiber connections being incident to nodes on the upper face $\nu \in \Nd_{B_u}$, i.e.,
\begin{subequations}
\begin{equation} \label{eq:test volume stress}
    T_{M}(\mathbf{x}) = -\sum\limits_{\nu \in \Nd_{B_u}} \sum\limits_{\mu \in \E(\nu)} \sum\limits_{r \in \F(\mu)} \mathbf{f}_{\mu, r}^{\nu}(\mathbf{x}) \cdot \mathbf{e_3}
\end{equation}
The strain $\epsilon_{M}$ of the structure corresponds to the change of height with regard to the initial structure height $h^\circ$. In the test it particularly reads, for $t\in (0,1]$
\begin{equation} \label{eq:test volume strain}
    \epsilon_{M}(t) = \frac{h(t) - h^\circ}{ h^\circ}= \frac{d}{h^\circ} t
\end{equation}
\end{subequations}
with the (maximal) displacement length $d$ as well as  $h(t) = (\max_{\nu \in \mathcal{N}_B} \mathbf{x}_{\nu}(t) - \min_{\nu \in \mathcal{N}_B} \mathbf{x}_{\nu}(t) ) \cdot \mathbf{e_3}$ and $h(0)=h^\circ$.
Consequently, the stress-strain relation is expressed by $(\epsilon_M(t), T_{M}( \mathbf{x}(\mathbf{z}(t),t))$, $t\in [0,1]$ where $\mathbf{z}(t)$ is the unique solution of the regularized quasi-static model. For an exemplary sample and \eqref{eq:dimensionless dynamic formulation}, the stress-strain relation with the respective material behavior during the tensile strength simulation is visualized in Fig.~\ref{fig:exemplary simulation}. The fiber structure shows the elongation with the expected simultaneous contraction.

\section{Graph-Based Data Reduction} \label{sec:Data Reduction}
The virtual graph-based fiber structures consist of several thousand nodes and edges. This results in a high dimensionality of the associated NESs and even leads to cases where the structural matrices cannot be stored. In this section, we present a strategy that accounts for the problem-relevant graph structure. It reduces data and, thus, the model's complexity and computing time, while achieving comparable tensile strength results. 

\begin{figure}[t]
\begin{subfigure}[c]{0.24\textwidth}
    \centering
    \includegraphics[width=\linewidth]{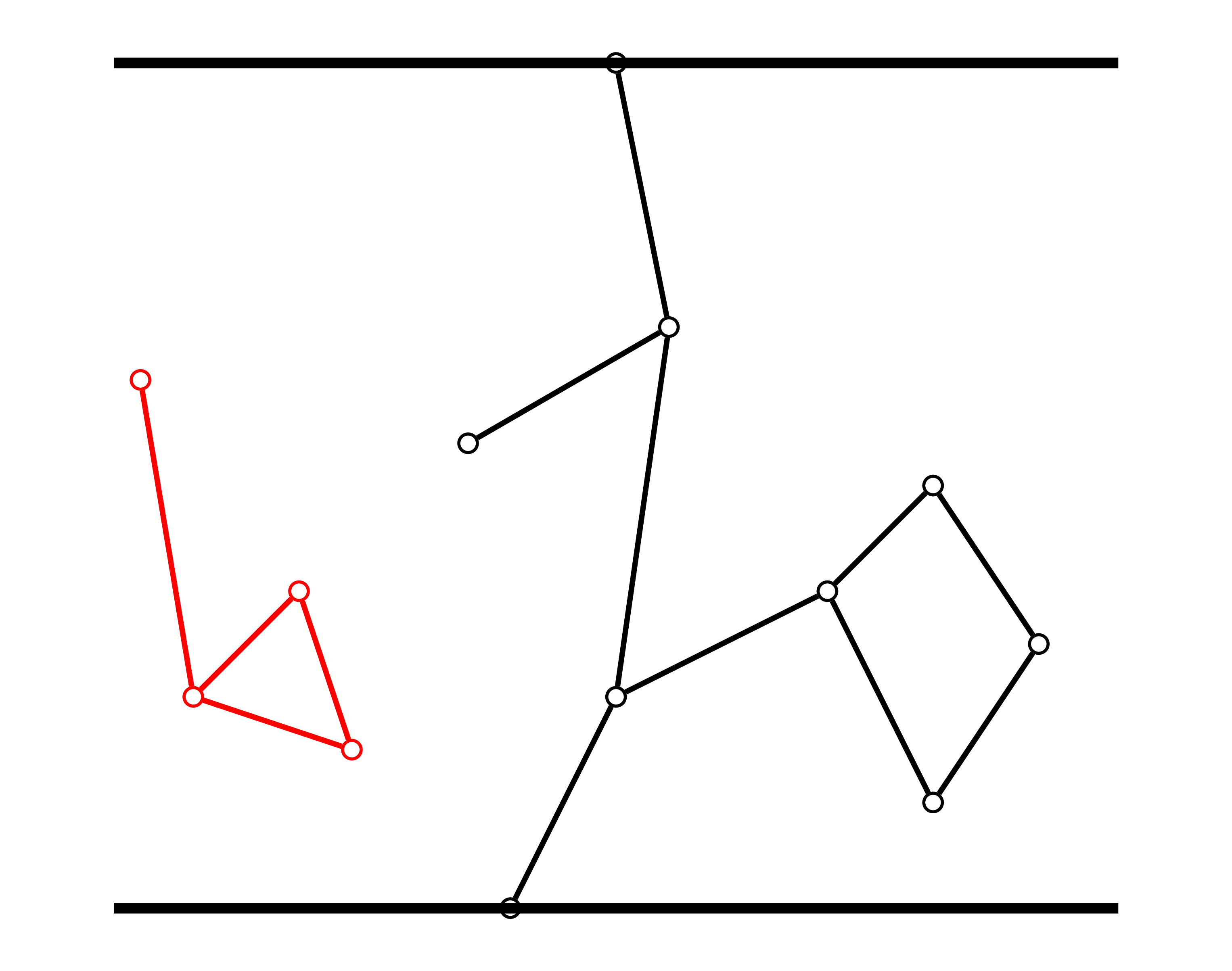}
    \caption{\tiny{Removal of uninvolved\\ \phantom{.......} components.}}
\end{subfigure}
\begin{subfigure}[c]{0.24\textwidth}
    \centering
    \includegraphics[width=\linewidth]{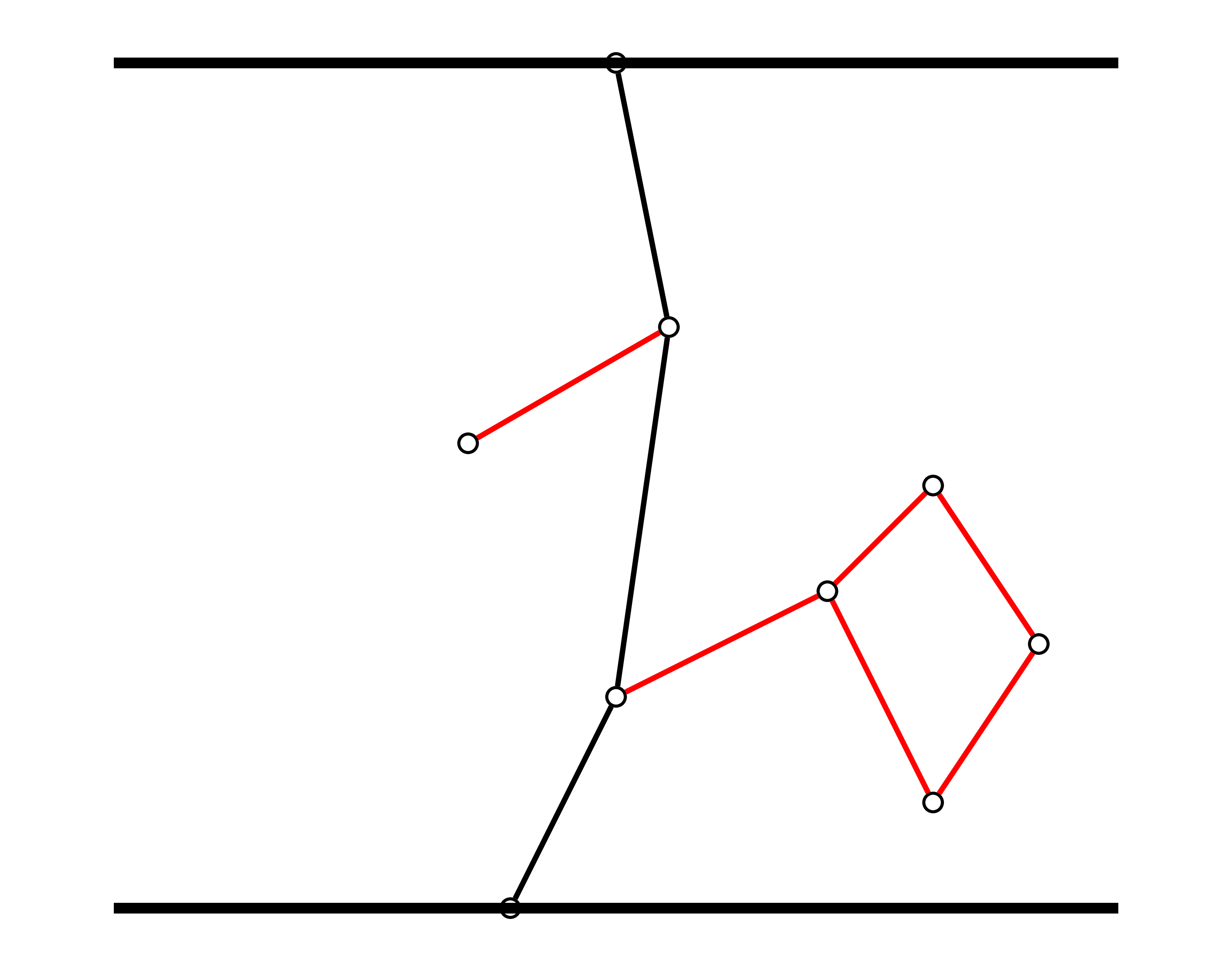}
    \caption{\tiny{Subsequent removal of\\ \phantom{.......} loose subgraphs.}}
\end{subfigure}
\begin{subfigure}[c]{0.24\textwidth}
    \centering
    \includegraphics[width=\linewidth]{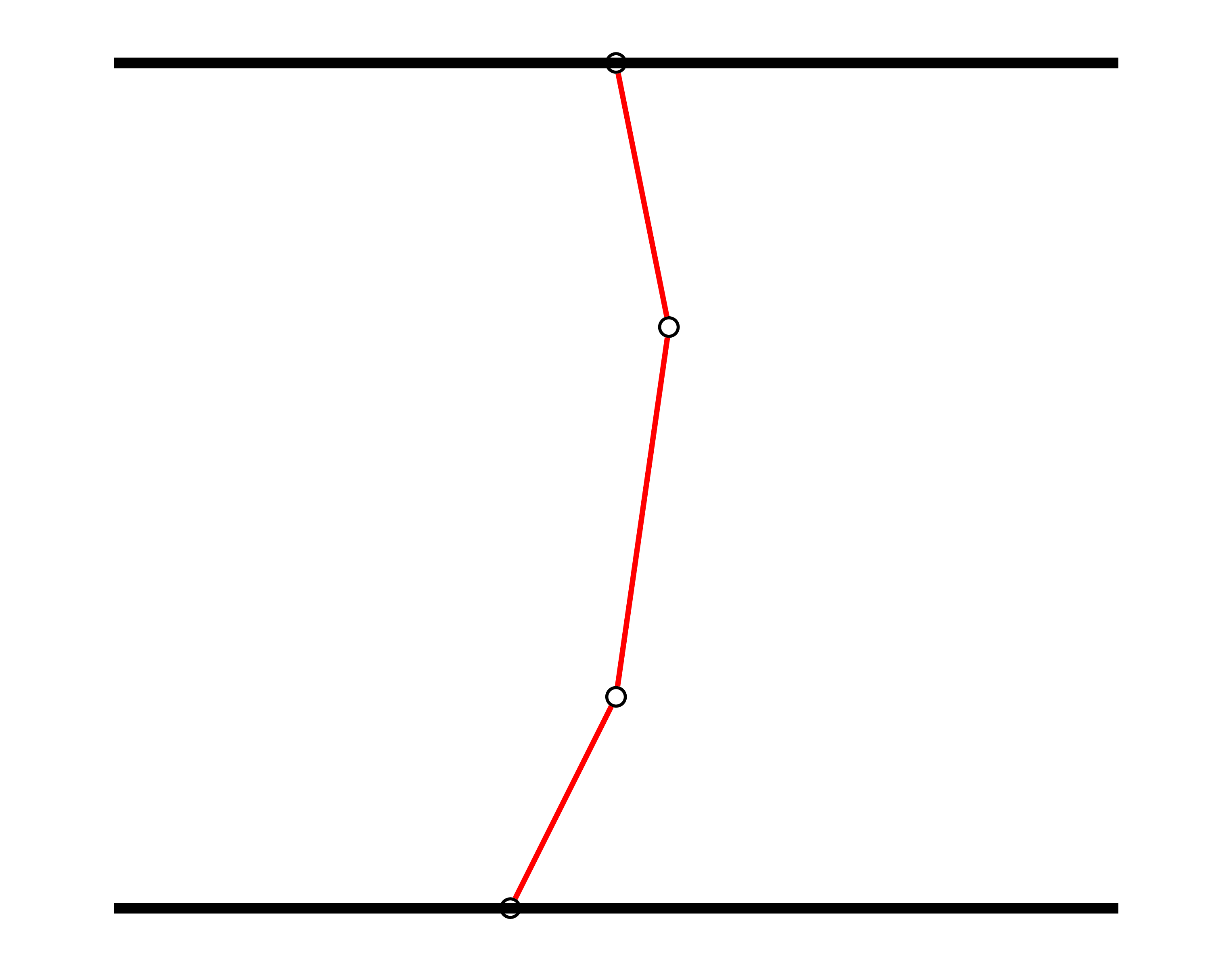}
    \caption{\tiny{Removal of simple linking\\ \phantom{.......} nodes, merging of fibers.}}
\end{subfigure}
\begin{subfigure}[c]{0.24\textwidth}
\centering
\includegraphics[width=\linewidth]{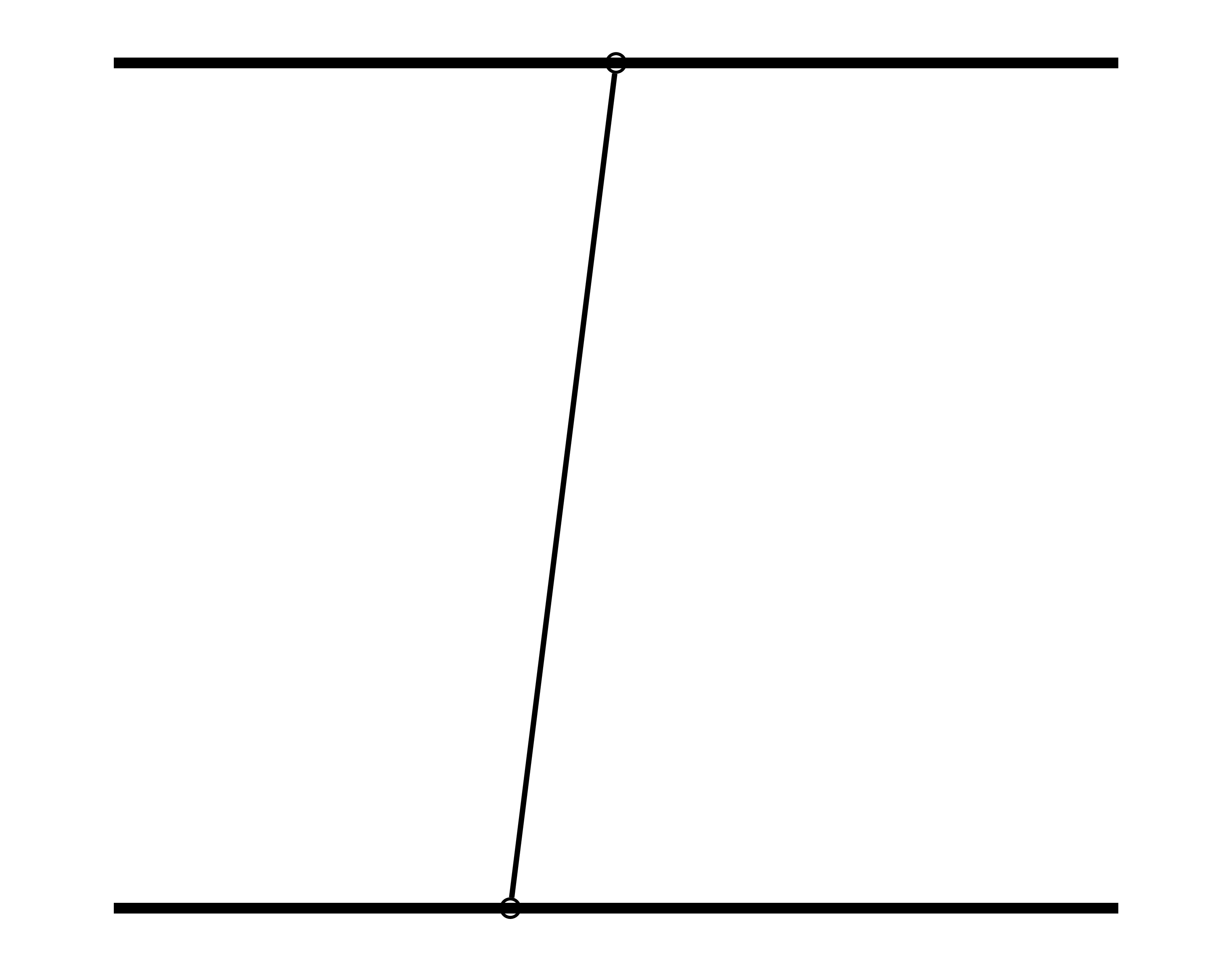}
\caption{\tiny{Resulting\\ \phantom{.......} reduced graph.}}
\end{subfigure}
\caption{Illustration of the graph-based data reduction from (A) to (D)}
\label{fig:data reduction strategies}
\end{figure}

Consider the underlying graph $\G$ to the fiber structure sample $M$. There exist subgraphs that do not contribute to the tensile strength behavior in the quasi-static (limit) model ($\varepsilon=0$), such as uninvolved components, loose subgraphs and simple linking nodes. Thus, we remove them from $\G$ in the subsequently described order. The tensile strength is determined by pulling away the upper and lower faces of the fiber structure. Concerning this procedure we distinguish between the involved and uninvolved graph components.  We refer to a graph component as uninvolved, $\mathcal{C}_U\subset \G$ , if it does not contain a path connecting both faces, i.e., $\mathcal{C}_U \cap \Nd_{B_u} = \emptyset$ or $\mathcal{C}_U \cap \Nd_{B_l} = \emptyset$. For each involved component $\mathcal{C}_I \subseteq \G$, a connected subgraph $\mathcal{S}\subset \mathcal{C}_I$ is said to be loose if it is connected to the remainder  $\mathcal{R} = \mathcal{C}_I \setminus \mathcal{S}$ over a cutvertex and if it does not contain a boundary node, $\Nd(\mathcal{S}) \cap \Nd_B = \emptyset$. Loose subgraphs do not give rise to any tension as they are pulled along over their respective cutvertices. Both, uninvolved components and loose subgraphs, can be found with a depth-first search \cite{cormen:b:2009}, which is coupled to an effort of $\mathcal{O}(|\Nd| + |\E|)$ comparisons. In addition, we remove all simple linking nodes in the involved graph components. We refer to a node $\nu$ as simple linking node, if it has degree 2 and links a fiber connection pair, i.e.,\ $|\E(\nu)|=2$ with $|\F(\mu)|=1$ for each $\mu\in \E(\nu)$. Replacing the original fiber connection pair by a single fiber connection of cumulated length yields a similar tensile strength behavior with less data. Simple linking nodes can be easily identified by means of the adjacency matrix and the edge-fiber connection relation, taking an effort of $\mathcal{O}(| \Nd |)$ comparisons. In Fig.~\ref{fig:data reduction strategies} the procedure of the data reduction is illustrated using a small exemplary graph.

In case of the regularized quasi-static model the fiber structure is exposed to friction ($\varepsilon >0$). We assume that the contributions of the uninvolved graph components, loose subgraphs and simple linking nodes to the overall tensile strength behavior are still negligibly small.  The removal of these subgraphs reduces crucially the size of the considered graph and the computing time and, hence, enables efficient tensile strength simulations. For a numerical study on the impact of the data reduction and the influence of $\varepsilon$ see Sec.~\ref{sec: numerics}.

\section{Numerical Results} \label{sec: numerics}

In this section we study the performance of the tensile strength simulations in view of accuracy and computational costs. Particular attention is paid to the influence of the regularizations and the data reduction. In addition, we investigate the effect of the randomness.

We consider a nonwoven material which is virtually generated with regard to the industrial scenario outlined in \cite{gramsch:p:2016}. Modeling the ramp-like contour (lay-down distribution) by a truncated normal distribution, the individual fibers in \eqref{sde:Final Model} are computed by means of the Euler-Maruyama scheme with step size $\Delta s$. For the virtual bonding \eqref{eq:contact point minimum}, we treat only a part of all fibers to be adhesive and to be able to form adhesive joints in order to reflect different fiber types in the production process. The contact threshold satisfies $\kappa < \Delta s$. The models and parameters used in the fiber structure generation are listed in Table~\ref{tab:parameter table}. In the tensile strength test we apply a maximal strain of $\epsilon^\star_M=0.5$ to the random graph-based fiber structure $\mathcal{M} = \{ \mathcal{G}, \F, \boldsymbol{\ell}, \mathbf{x} \}$. To account for the stiffness of the regularized quasi-static model \eqref{eq:dimensionless dynamic formulation} we employ a BDF(2)-scheme (multistep method, Backward Differentiation Formula) with step size $\Delta t$ which is initialized with an implicit midpoint rule. The resulting nonlinear equation systems are solved with an exact Newton method with an error tolerance of $\mathcal{O}(10^{-8})$.  As initial guess we use a linear extrapolation of the variables, i.e., ${\mathbf{z}}^{(0)}_{k+1} = 2 \mathbf{z}_k - \mathbf{z}_{k-1}$, as warm start. Note that the involved linear systems are sparse due to the modeled material law and tackled with a direct solver using analytically determined Jacobians.  All computations have been carried out with MATLAB version R2019b on a Fujitsu Esprimo P958 running Ubuntu 20.04.

\begin{table}[tb]
\begin{tabular}{llll}
\toprule
\bf Property & \bf Symbol &\multicolumn{2}{l}{\bf Model}     \\
\midrule
Fiber lay-down distribution & $r $ &\multicolumn{2}{l}{$r(x)=\exp(-\frac{x^2}{2 \sigma^2}) \chi_{[-\frac{x_r}{2}, \frac{x_r}{2}]}/\int_{-\frac{x_r}{2}}^{\frac{x_r}{2}} \exp(-\frac{x^2}{2 \sigma^2}) \mathrm{d}x$}\\
Production time for $\mathcal{V}_R$ & $T_R$  &\multicolumn{2}{l}{$T_R=(x_r+w_R)/v_B$}\\
Trace curve due to belt motion \hspace*{0.3cm}& $x_B$ &\multicolumn{2}{l}{$x_B(T)=(x_r-w_R)/2+v_BT$}\\
Lay-down potential & $V$ & \multicolumn{2}{l}{$V(\boldsymbol{\xi})=\boldsymbol{\xi} \cdot (\sum_{i=x,y,z} (\sigma_i^2)^{-1} \mathbf{e_i} \otimes \mathbf{e_i}) \cdot \boldsymbol{\xi}$}\\
\midrule
\midrule
&   &  {\bf Value} & {\bf Unit}  \\
\midrule
Test volume height & $h$ & $5.0 \cdot 10^{-2} $ &m\\
Test volume width & $w, w_R$ &  $1.0 \cdot 10^{-2}, 1.2 \cdot 10^{-1}$ &m\\
Fiber length  &  $L$ &  $5.5 \cdot 10^{-2}$ & m  \\
Fiber number (adhesive part) for $\mathcal{V}_R$& $n$ & $7.3 \cdot 10^{5} \quad (40\%)$ &\\
Lay-down distribution deviation  & $\sigma$ & $ 2.0\cdot 10^{-2}$ & m\\
Lay-down distribution range & $ x_r$ & $2.0\cdot 10^{-1}$ &m\\
Conveyor belt speed & $v_B$ & $3.3\cdot 10^{-2}$ & m\,s$^{-1}$\\
Noise amplitude &  $A$ & $2.0 \cdot 10^{-1}$ & m$^{-1/2}$\\
Anisotropy constant &  $B$ & $3.0 \cdot 10^{-1}$ &\\
Lay-down potential parameters & $\sigma_x,\sigma_y , \sigma_z$  & $1.5 \cdot 10^{-2}, 2.0 \cdot 10^{-2}, 1.5 \cdot 10^{-3}$ &m\\
Contact threshold & $\kappa$ & $2.8 \cdot 10^{-4}$ & m\\
Fiber discretization length & $ \Delta s$ & $ 3.7 \cdot 10^{-4} $ & m\\
Elasticity modulus &$ EA$& $1.0$ & N\\
Maximal displacement length & $d$ & $2.5\cdot 10^{-2}$ & m\\
\bottomrule
\end{tabular} 
\caption{Models and parameters used for virtual fiber structure generation and tensile strength test. } 
\label{tab:parameter table}
\end{table}
  
\subsection{Impact of regularizations} \label{sec:numerical results}
A solution of the quasi-static model and hence the tensile strength behavior depends on the choice of the two involved regularization parameters:  $\varepsilon$ (friction) in \eqref{eq:dimensionless dynamic formulation} and $\delta$ (material law) in \eqref{eq:stress-strain relation}. To study their impact on the solution and to exclude randomness we consider a fixed (deterministic) material sample $M=\mathcal{M}(\omega)$ and $\Delta t=10^{-5}$ in the tensile strength test.  

\begin{figure}
\begin{subfigure}[t]{0.30\textwidth}
    \centering
    \includegraphics[width=\linewidth]{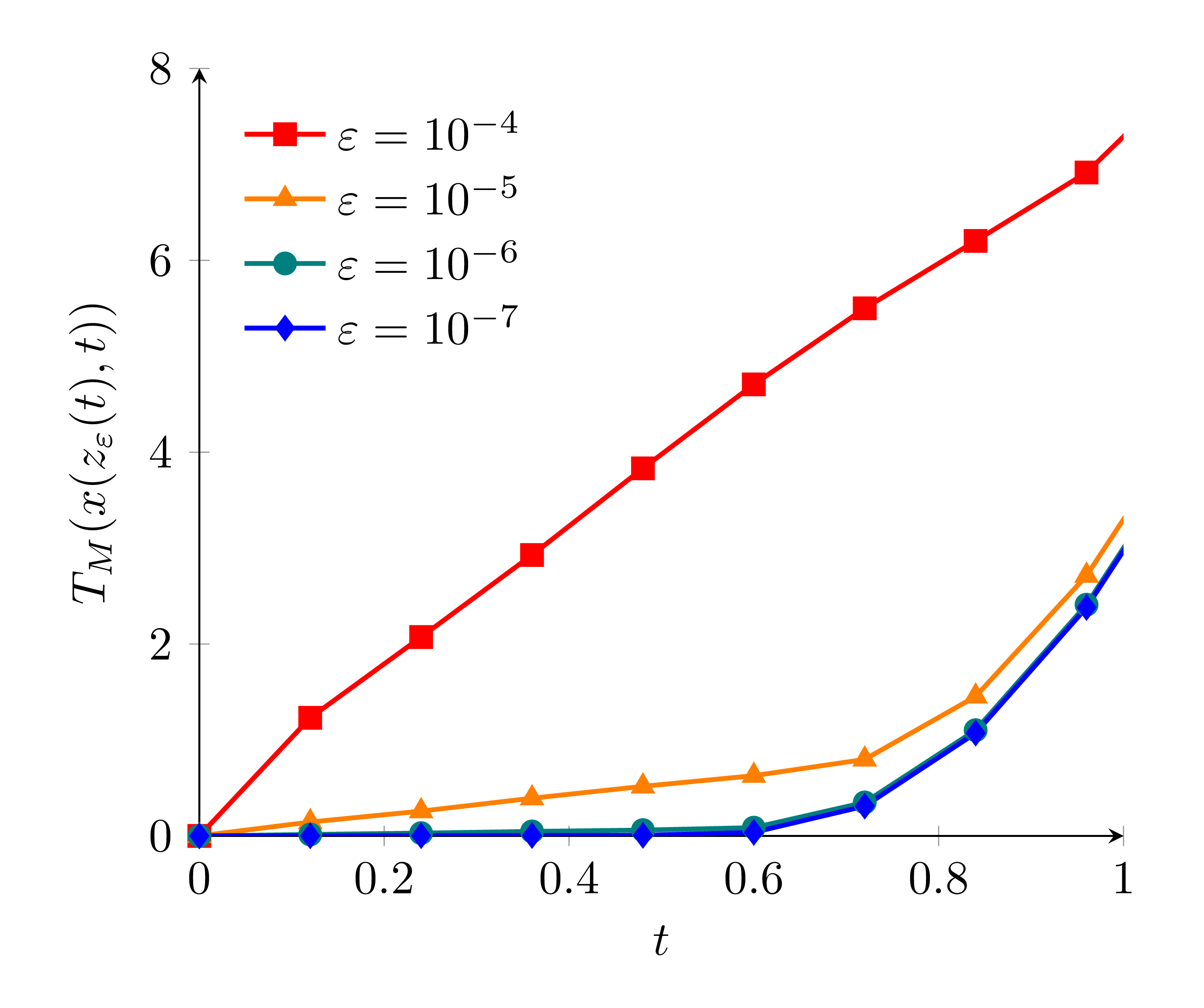}
    \caption{\tiny Tensile force curve corresponding\\ \phantom{.......} to the solution $\mathbf{z}_{\varepsilon}$ of \eqref{eq:dimensionless dynamic formulation}.}
    \label{fig:Convergence plot}
\end{subfigure}
\hfill
\begin{subfigure}[t]{0.32\textwidth}
    \centering
    \includegraphics[width=\linewidth]{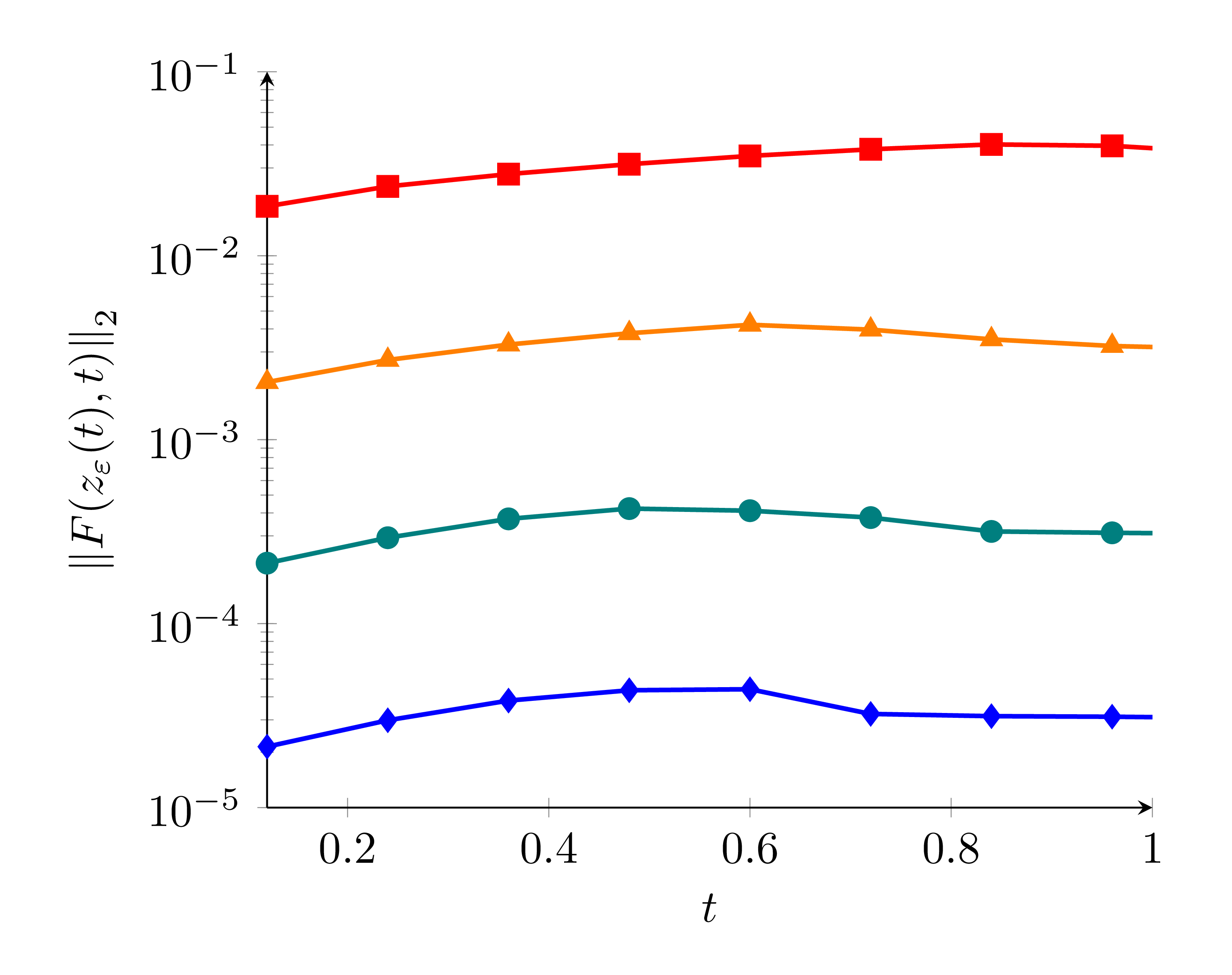}
    \caption{\tiny Residual $\|\mathbf{F}(\mathbf{z}_\varepsilon(t),t)\|_2$ for various $\varepsilon$.\\ \phantom{.......} }
    \label{fig:quasi-static model error}
\end{subfigure}
\hfill
\begin{subfigure}[t]{0.32\textwidth}
    \centering
    \includegraphics[width=\linewidth]{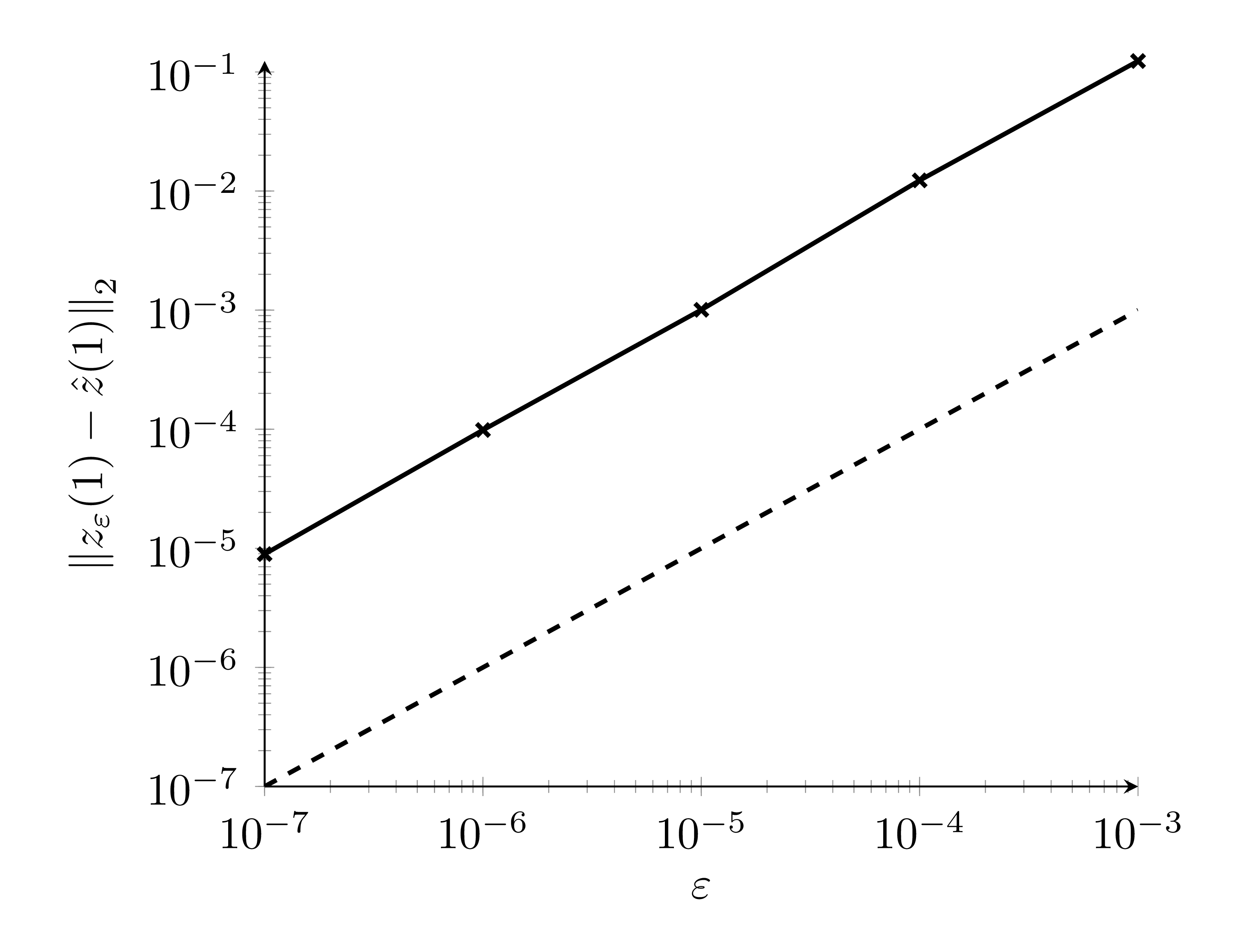}
    \caption{\tiny Error to reference solution $\hat{\mathbf{z}}(1)$. Dash-\\ \phantom{.......} ed line shows first order convergence.}
    \label{fig:regularization effect}
\end{subfigure}
\caption{\protect\centering Influence of $\varepsilon$ on the tensile strength behavior, $\delta=10^{-4}$.}
\label{fig:impact_eps}
\end{figure}

Figure~\ref{fig:impact_eps} illustrates the influence of the friction-associated regularization parameter $\varepsilon$ for fixed $\delta=10^{-4}$, where the reference solution $\hat{\mathbf{z}}$ is computed with $\varepsilon=10^{-8}$. The visualized curves for the tensile forces over $t\in [0,1]$ can be viewed as scaled stress-strain relations, they correspond to solutions $\mathbf{z}_{\varepsilon}$ of \eqref{eq:dimensionless dynamic formulation} for $\varepsilon \in \{ 10^{-4},10^{-5},10^{-6}, 10^{-7} \}$. For $\varepsilon \to 0$,   $\mathbf{z}_{\varepsilon}$ converges to the friction-associated solution of the quasi-static limit model $\mathbf{F}(\mathbf{z},t)=\mathbf{0}$. The convergence turns out to be linear in $\epsilon$. Over the tensile strength test $t\in [0,1]$, the residuals of NES show an almost constant magnitude of $\mathcal{O}(100\,\varepsilon)$ for given $\varepsilon$. For moderate $\varepsilon$, $\varepsilon\geq 10^{-4}$, this implies large absolute and relative errors in the beginning of the test (i.e., for small $t$). The errors accumulate for increasing $t$ and lead to undesirably high deviations from the reference in the final node constellation at $t=1$, e.g., $ \|\mathbf{z}_\varepsilon(1) - \hat{\mathbf{z}}(1)\|_2 \sim \mathcal{O}(10^1)$ for $\varepsilon= 10^{-4}$. Moreover, the corresponding stress-strain curves show a qualitatively absolutely wrong behavior. However, the regularization parameter $\varepsilon$ cannot be chosen arbitrarily small either, since due to the stiffness of the model the computational costs grow significantly for $\varepsilon \to 0$. A good compromise between accuracy and effort turns out to be $\varepsilon = 10^{-6}$, because it provides satisfactory approximations with an acceptable computational effort as we will comment on later on.

\begin{figure}
    \begin{subfigure}[t]{0.32\textwidth}
        \centering
        \includegraphics[width=\linewidth]{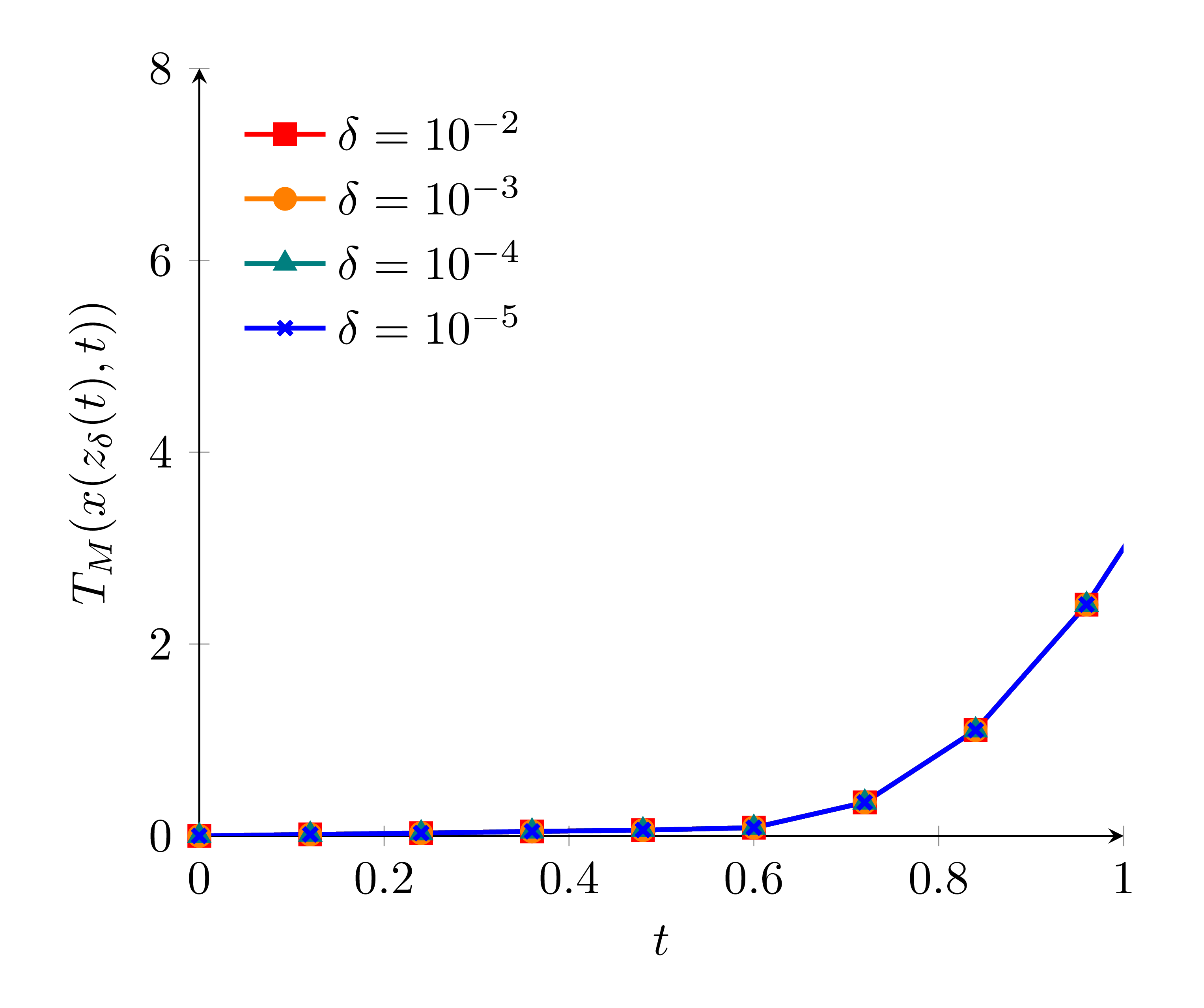}
        \caption{\tiny Tensile force curve corresponding\\ \phantom{.......} to the solution $\mathbf{z}_{\delta}$ of \eqref{eq:dimensionless dynamic formulation}.}
        \label{fig:delta stress-strain curves}
    \end{subfigure}
    \hfill
    \begin{subfigure}[t]{0.34\textwidth}
        \centering
        \includegraphics[width=\linewidth]{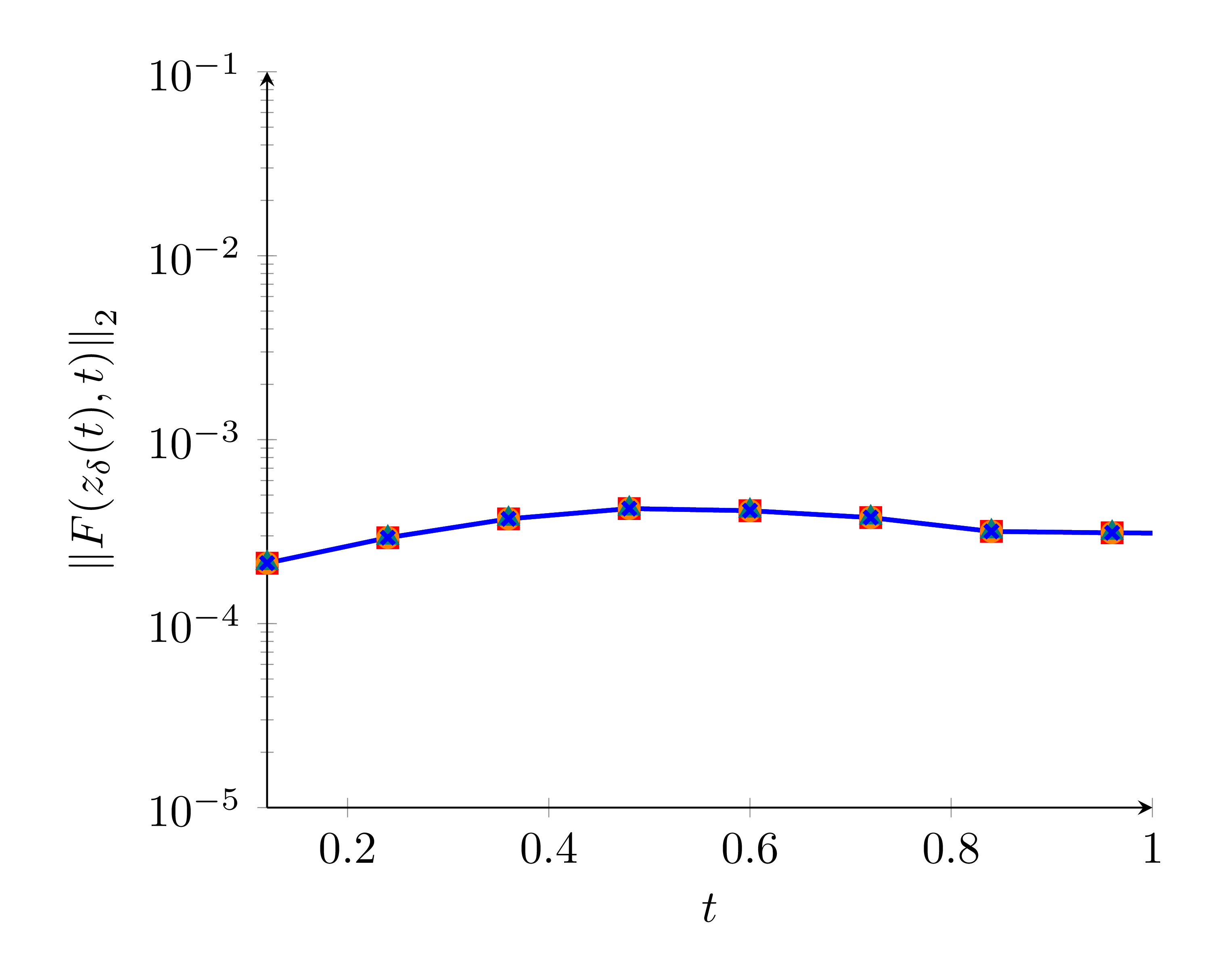}
        \caption{\tiny Residual $\|\mathbf{F}(\mathbf{z}_\delta(t),t)\|_2$ for various $\delta$.\\ \phantom{.......}}
        \label{fig:delta error to quasi-static model}
    \end{subfigure}
    \hfill
    \begin{subfigure}[t]{0.32\textwidth}
        \centering
        \includegraphics[width=\linewidth]{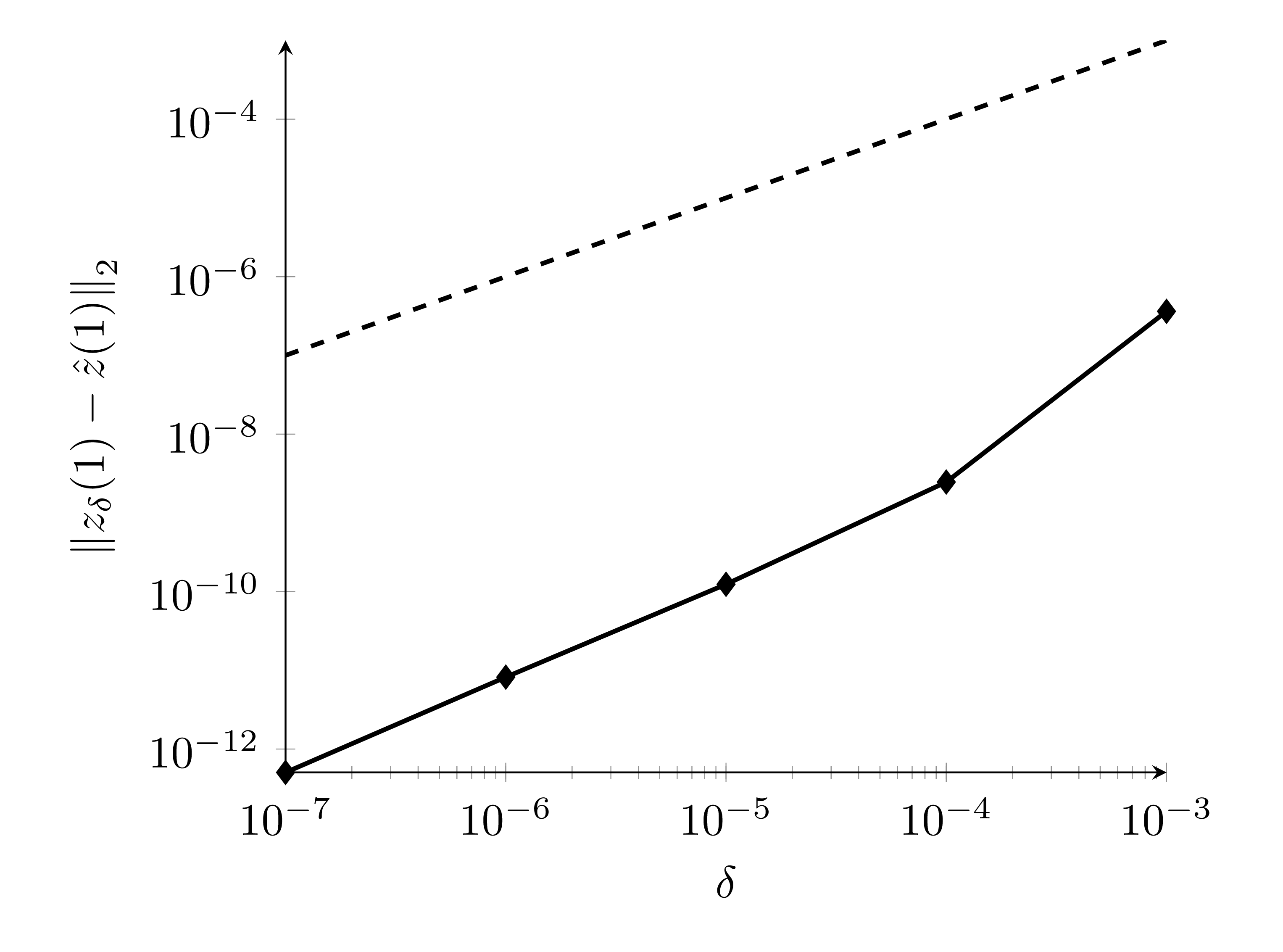}
        \caption{\tiny Error to reference solution $\hat{\mathbf{z}}(1)$. Dash-\\ \phantom{.......} ed line shows first order convergence.}
        \label{fig:delta error behavior}
    \end{subfigure}
    \caption{\protect\centering Influence of $\delta$ on the tensile strength behavior, $\varepsilon=10^{-6}$.}
    \label{fig:impact_delta}
\end{figure}

The parameter $\delta$ acts as smoothing parameter in the underlying material law. Thus, it is not only included in the regularized quasi-static model \eqref{eq:dimensionless dynamic formulation}, but also in the limit model $\mathbf{F}(\mathbf{z},t)=\mathbf{0}$ itself. Figure~\ref{fig:impact_delta} illustrates the influence of $\delta$ for fixed $\varepsilon=10^{-6}$, where the reference solution $\hat{\mathbf{z}}$ is computed with $\delta=10^{-8}$. 
The visualized tensile force curves correspond to solutions $\mathbf{z}_{\delta}$ of \eqref{eq:dimensionless dynamic formulation} for $\delta \in \{10^{-2},10^{-3},10^{-4},10^{-5}\}$. Apparently, the sensitivity of the stress-strain results towards a variation of $\delta$ is comparably small. The convergence to the reference is linear in $\delta$, satisfactory results are already obtained for moderate $\delta$, e.g., $ \|\mathbf{z}_\delta(1) - \hat{\mathbf{z}}(1)\|_2 \sim \mathcal{O}(10^{-6})$ for $\delta= 10^{-4}$. The residuals of NES show the same small constant magnitude over $t$ independent of $\delta$.
Thus, the impact of $\delta$ can be neglected, we use $\delta = 10^{-4}$ in the following.

\subsection{Effects of randomness} \label{sec:randomness in the simulations}
To investigate the effects of randomness on the tensile strength behavior of the fiber structure $\mathcal{M}$, 
we perform Monte-Carlo simulations. The focus is on the performance of the data reduction in view of accuracy and computational savings. The presented results are based on $n_\mathcal{M}$ samples, using $n_\mathcal{M}=100$, $\varepsilon=10^{-6}$, $\delta=10^{-4}$ and $\Delta t=10^{-4}$.

\begin{figure}
\begin{subfigure}[t]{0.32\textwidth}
    \centering
    \includegraphics[width=\linewidth]{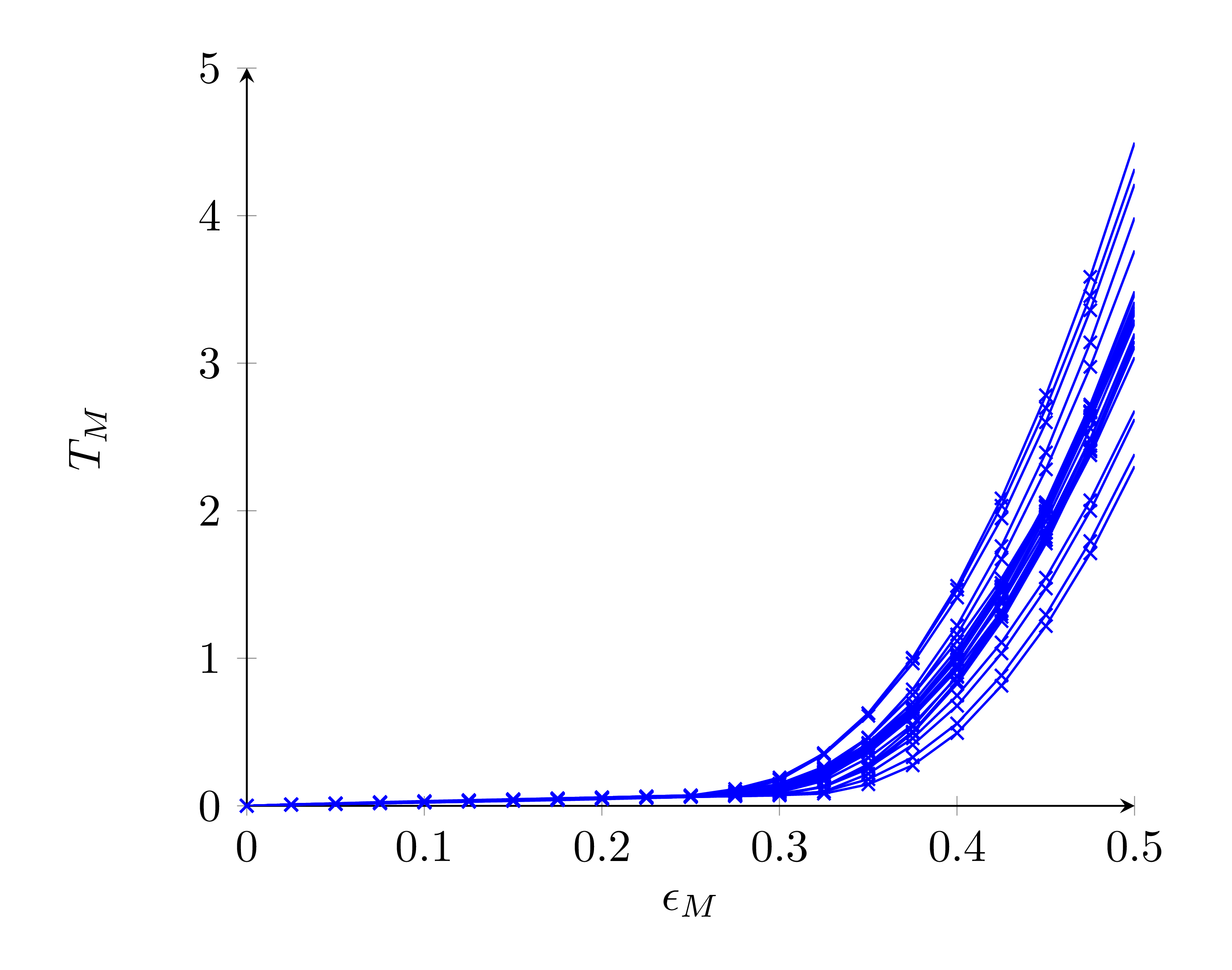}
    \caption{\tiny Stress-strain relation.}
    \label{fig:mc stress-strain curves}
\end{subfigure}
\hfill
\begin{subfigure}[t]{0.32\textwidth}
    \centering
    \includegraphics[width=\linewidth]{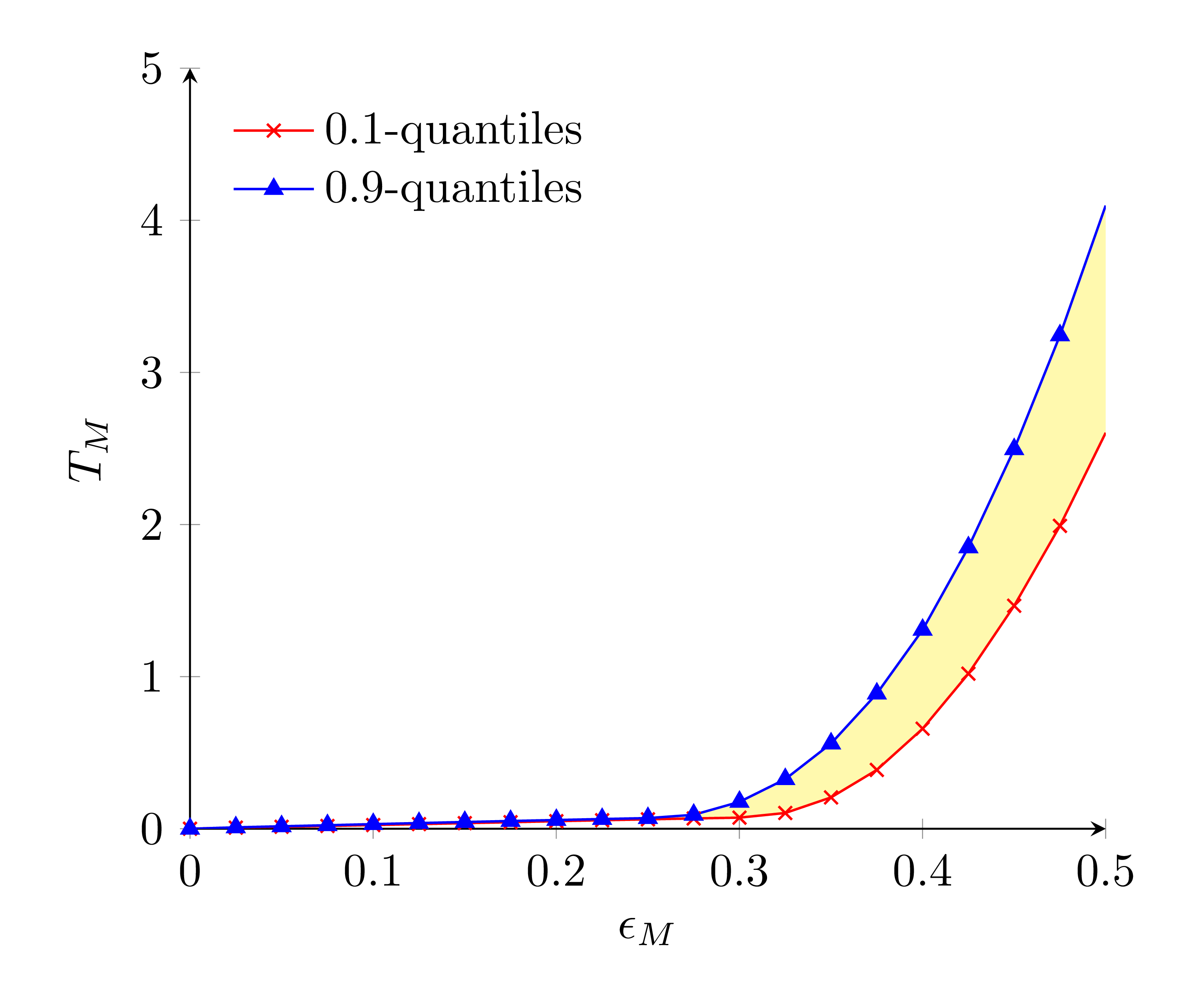}
    \caption{\tiny $0.1$-quantile and $0.9$-quantile.}
    \label{fig:mc quantiles}
\end{subfigure}
\hfill
\begin{subfigure}[t]{0.32\textwidth}
    \centering
    \includegraphics[width=\linewidth]{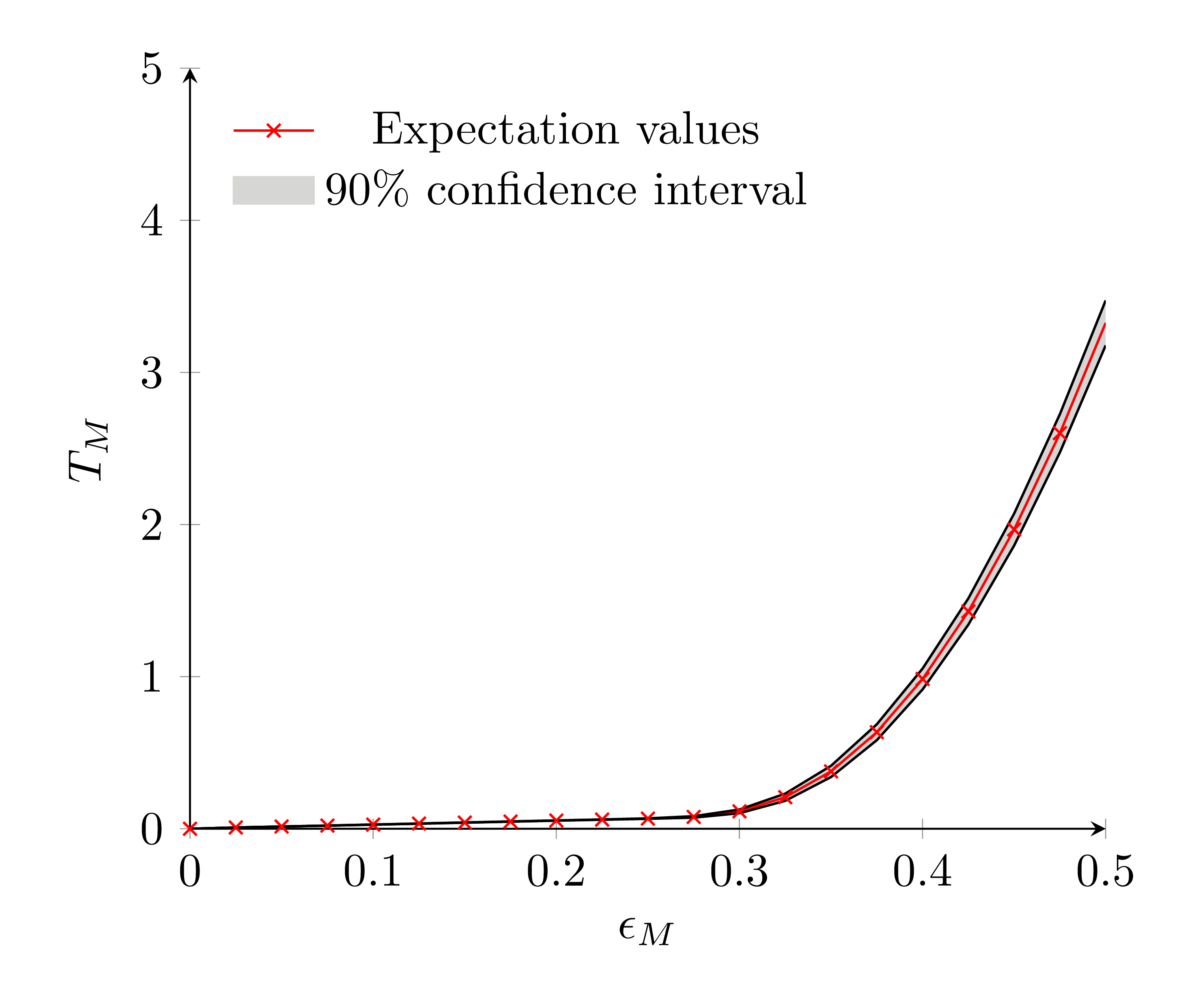}
    \caption{\tiny Average and $90\%$ confidence interval.}
    \label{fig:mc expectated behavior}
\end{subfigure}
\caption{Monte-Carlo results for the stress-strain relation of the (full) fiber structure $\mathcal{M}$, $n_\mathcal{M}=100$.}
\label{fig:monte-carlo experiment}
\end{figure}

The Monte-Carlo results for the stress-strain relations $(\epsilon_M(t),T_M(\mathbf{x}(\mathbf{z}(t),t))$, $t\in[0,1]$, $M=1,...,n_\mathcal{M}$ of the random fiber structure are visualized in Fig.~\ref{fig:monte-carlo experiment}. The results show the features that are known and expected from experimental data. To assess characteristic material properties it is convenient to consider quantiles and average. The computation of the confidence interval of the expected tensile strength behavior is based on the assumption that the results are normally distributed. Apart from a few outliers the stress-strain curves are very similar which leads to tight confidence intervals as desired in the industrial production process. The results refer to samples without data reduction.

\begin{figure}
\begin{subfigure}[t]{0.32\textwidth}
    \centering
    \includegraphics[width=\linewidth]{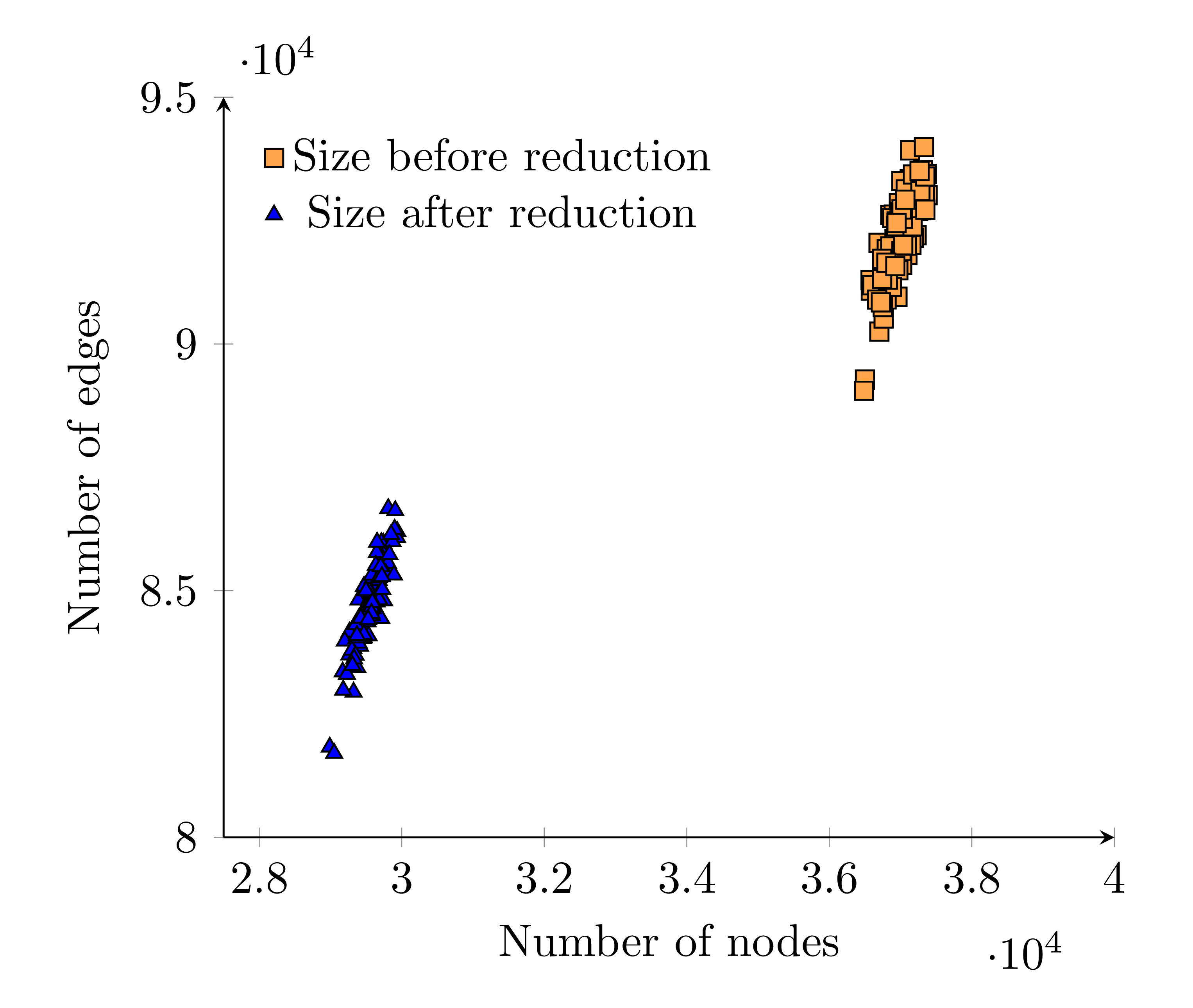}
    \caption{\tiny Graph size before and after reduction.}
    \label{fig:size reduction}
\end{subfigure}
\hfill
\begin{subfigure}[t]{0.32\textwidth}
    \centering
    \includegraphics[width=\linewidth]{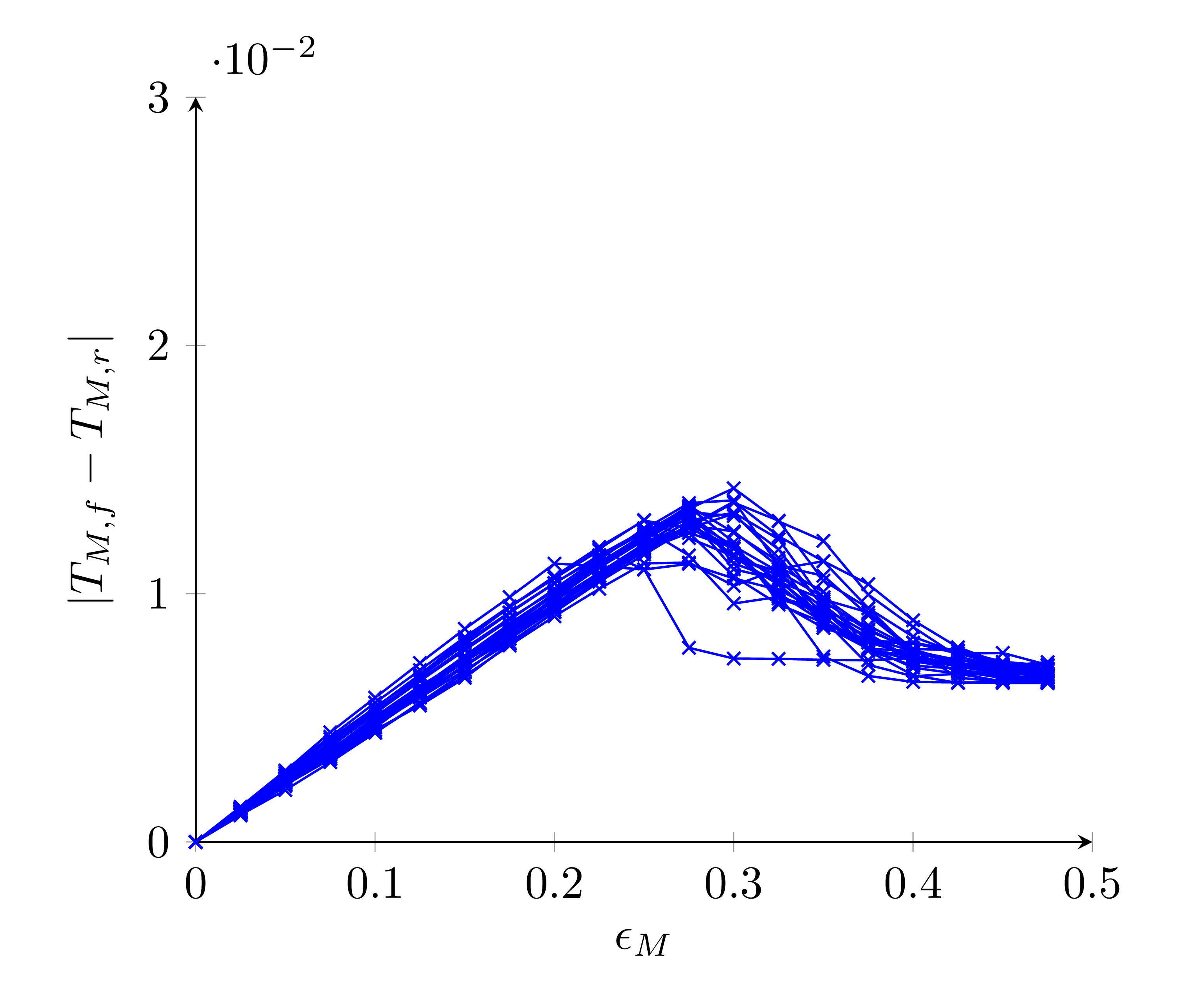}
    \caption{\tiny Absolute reduction error in tensile\\ \phantom{.......} forces.}
    \label{fig:reduction errors}
\end{subfigure}
\hfill
\begin{subfigure}[t]{0.32\textwidth}
    \centering
    \includegraphics[width=\linewidth]{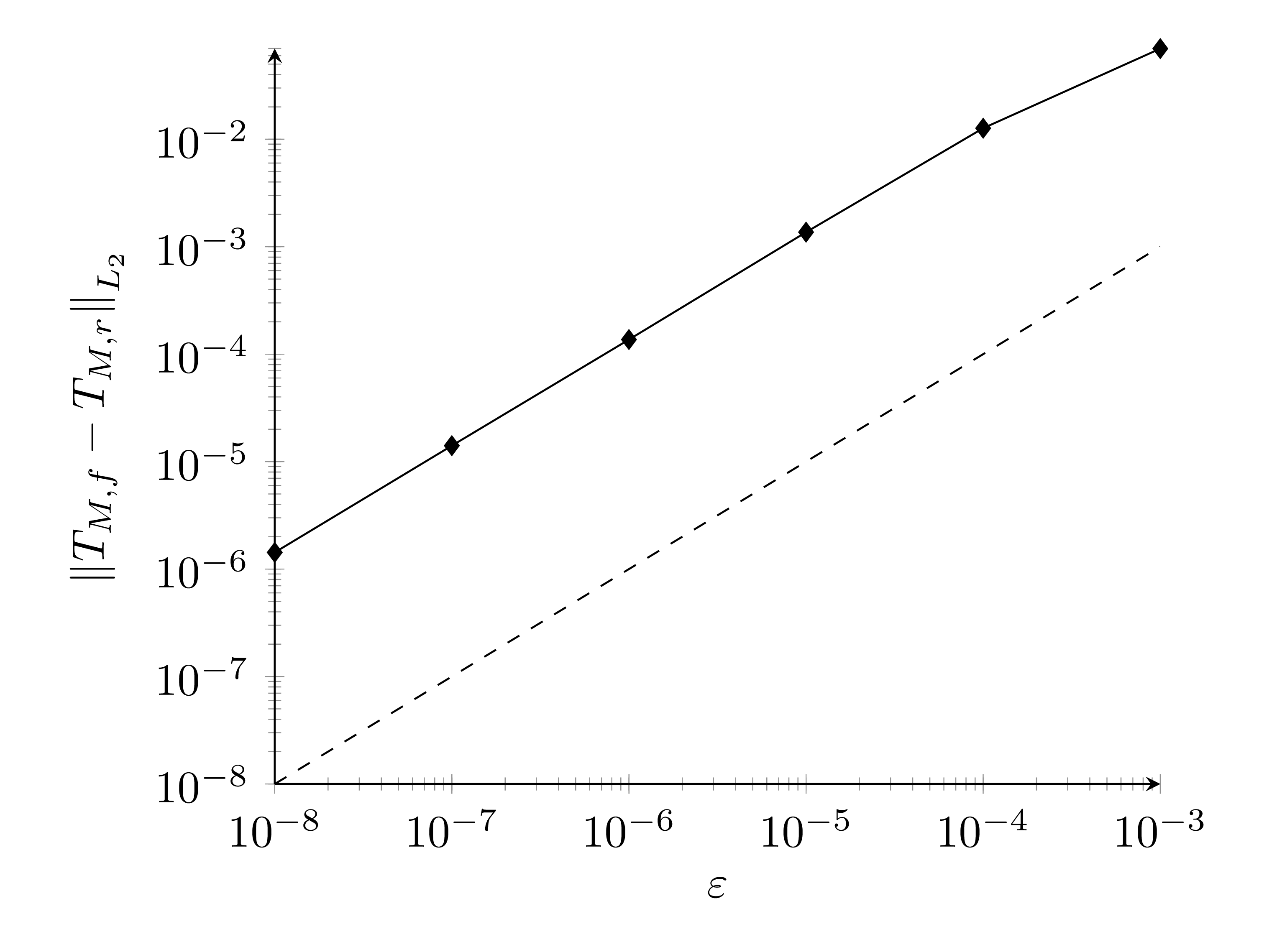}
    \caption{\tiny $L^2$-reduction error depending on $\varepsilon$.\\ \phantom{......} First order convergence by dashed line.}
    \label{fig:vanishing reduction error}
\end{subfigure}
\caption{\protect\centering Effect of data reduction on the tensile strength behavior, cf.~Fig.~\ref{fig:monte-carlo experiment}.}
\label{fig:data reduction influence}
\end{figure}

Depending on the random structure of the underlying graphs, the amount of nodes and edges that is removed by our data reduction procedure varies for each sample.  In the considered setup the (full) graph has in average $3.8\cdot 10^{4}$ nodes and $9.3 \cdot 10^{4}$ edges, the average reduction of nodes is approximately $20\%$ and that of edges $8~\%$, see Fig.~\ref{fig:size reduction} for details on the distribution. Thus, by data reduction the dimension of the problem is significantly decreased. Comparing the tensile strength behavior before and after data reduction, we refer to the tensile force of the full sample as $T_{M,f}$ and to the one of the reduced sample as $T_{M,r}$. The reduction errors throughout the tensile strength test are of moderate, acceptable magnitude in total (i.e., $L^2$-error of size $\mathcal{O}(10^{-4})$ for $\varepsilon=10^{-6}$). The highest absolute deviation is observed in the stress-strain situation when the tensile force starts to grow rapidly for increasing strain $\epsilon_M$. However, note that the relative deviations are large in the beginning of the tensile strength test when the fiber structure is in an (almost) stress-free state.  Reason is that -- by removing uninvolved components, loose subgraphs and simply linked nodes -- the frictional effects on the respective nodes due to the introduced regularization parameter $\varepsilon$ are neglected. For $\varepsilon \to 0$ these effects vanish, as the removed parts do not give rise to tension anymore. In particular, the $L^2$-reduction error goes linearly in $\varepsilon$, cf.~Fig.~\ref{fig:data reduction influence}.

Typical computation times recorded for sample generation, data reduction and the corresponding tensile strength simulations can be found in Table~\ref{tab:efficiency}. Apparently, the fiber structure generation, which is mainly determined by the identification of adhesive joints, is comparably fast. In this context, the effort for data reduction is marginal and can be neglected. Most of the time is spent for the actual tensile strength simulations.  Solving the emerging nonlinear systems is, hereby, the computational bottleneck. To ensure (numerical) convergence of Newton's method a good initial guess is necessary which requires here in general a small step size $\Delta t$. The interplay of step size and number of Newton iterations per step determines the overall computation time. We employ step size control to rise the efficiency. The step size is chosen in the range $\Delta t^\star \in[10^{-6},10^{-2}]$ to satisfy an accuracy $\mathcal{O}(10^{-8})$. The use of parallelization is beneficial for the Monte-Carlo simulations.

\begin{table} \centering
\begin{tabular}{l|c|ccc|ccc}
\toprule
{\bf Procedure}  & {\bf DOF} & \multicolumn{3}{c|}{{\bf Newton iterations per step}$^\diamond$} & \multicolumn{3}{c} {\bf CPU [min]}  \\
$\Delta t$ & & $10^{-4}$ & $10^{-5}$ & $\Delta t^\star$ & $10^{-4}$ & $10^{-5}$ & $\Delta t^\star$\\
\midrule
structure generation &&&&&  $\phantom{0}20$ & $\phantom{0}20$ & $\phantom{0}20$ \\
full simulation & 100\,000   & 8.2 & 3.2 & 10.0 & 340 & 1280 & 130\\
data reduction &&&&&  $\phantom{00}2$ & $\phantom{00}2$ & \phantom{00}2 \\
reduced simulation & \phantom{0}80\,000 & 8.2 & 3.1 & 11.5 & 280 & 1080 & 100 \\ 
\bottomrule
\end{tabular} 
\caption{\protect\centering Typical computation times for a sample, ($^\diamond$ averaged value).} 
\label{tab:efficiency}
\end{table}

This proof of concept shows that tensile strength simulations of nonwoven materials can be performed at a macroscopic level without the use of additional homogenization techniques. This way, structural features coming from the production process (like the ramp-contour) are taken into account.

\section{Conclusion and Outlook}

In this paper we presented a framework that allows for efficient virtual tensile strength investigations on nonwoven materials in which structural production-related features are taken into account on a macroscopic level without any homogenization techniques. Modeling the large complex elastic random fiber structure as graph-based and of truss-type, we introduced a friction-associated regularization and developed a problem-tailored data reduction method. The randomness is accounted for by Monte-Carlo simulations. We performed a proof of concept for an industrial setup. The computed stress-strain curves show the typical properties known from measurements. For comparability with experiments and validation issues, however, there is need of further extensive and careful calibrations of the virtual fiber structure generation with tomographic data.
A longtime objective is the embedding of the simulations in an optimization framework in order to control the production parameters to obtain a desired material behavior with respect to minimal material usage, or to minimize the variance of the stress-strain curves.

\subsection*{Acknowledgements} The support of the German Research Training Group \emph{Algorithmic Optimization -- ALOP} (RTG 2126) funded by the DFG is acknowledged.

\bibliographystyle{siam}
\bibliography{paper}

\end{document}